\documentclass[11pt,reqno]{amsart}
\usepackage{fixltx2e}                    
\usepackage[osf]{mathpazo}  
\usepackage[utf8]{inputenc}            
\usepackage{amsmath}                     
\usepackage{amssymb, latexsym, stmaryrd, amsthm, dsfont, amsfonts, amsbsy,amsthm, amsmath, mathrsfs}            
\usepackage{mathtools}                   
\usepackage{bm}                          
\usepackage{enumerate}                   
\usepackage{verbatim}                    
\usepackage{url}   
\usepackage{lscape}                      
\usepackage{soul}

\usepackage[margin=1.3in]{geometry}

\usepackage{microtype}  
\usepackage[all, knot]{xy}
        \xyoption{arc} 
        \xyoption{web}                 
\makeatletter                            
\def\MT@register@subst@font{\MT@exp@one@n\MT@in@clist\font@name\MT@font@list
 \ifMT@inlist@\else\xdef\MT@font@list{\MT@font@list\font@name,}\fi}
\makeatother

\usepackage[pdftex,bookmarks,bookmarksnumbered,linktocpage,   
         colorlinks,linkcolor=blue,citecolor=blue]{hyperref}

\newcommand{\bit}{\begin{itemize}}    
\newcommand{\eit}{\end{itemize}}
\newcommand{\ben}{\begin{enumerate}}
\newcommand{\een}{\end{enumerate}}

\newcommand{\benroman}{\ben[\normalfont (i)]}  
\let\eroman\een

\newcommand{\bde}{\begin{description}}
\newcommand{\ede}{\end{description}}

\newcommand{\?}{\ensuremath{\mkern0.4\thinmuskip}}   

\let\leq=\leqslant
\let\geq=\geqslant

\let\epsilon=\varepsilon
\let\Lambda\varLambda
\let\Gamma\varGamma
\let\Delta\varDelta
\let\Lambda\varLambda
\let\Omega\varOmega
\let\Theta\varTheta
\let\Xi\varXi
\let\Pi\varPi
\let\Sigma\varSigma

\let\class=\mathsf                              
\let\oper=\mathbb                               

\bmdefine{\A}{A}                                
\bmdefine{\2}{2}
\bmdefine{\B}{B}
\bmdefine{\D}{D}
\bmdefine{\M}{M}                                
\bmdefine{\LLL}{L}                              
\bmdefine{\Fm}{Fm}                              
\bmdefine{\zerou}{[0{,}1]}  
\bmdefine{\T}{T}                                

\newcommand{\LL}{\mathscr{L}}                   

\newcommand{\PPP}{\oper{P}}
\newcommand{\PSD}{\oper{P}_{\!\textsc{sd}}^{}}

\newcommand{\PRk}{\oper{P}_{\!\textsc{r}_{\kappa^{+}}}^{}}

\newcommand{\SSS}{\oper{S}}
\newcommand{\III}{\oper{I}}
\newcommand{\UUU}{\oper{U}}
\newcommand{\RRR}{\oper{R}}

\newcommand{\Con}{\mathrm{Con}}                            
      
\newcommand{\Mod}{\class{Mod}}

\bmdefine{\boldstar}{\mathchoice{\textstyle*}{\textstyle*}{\textstyle*}{\scriptstyle*}}

\newcommand{\ModS}{\Mod^{\equiv}}

\bmdefine{\btau}{\tau}                                  
\bmdefine{\brho}{\rho}                                  

\bmdefine{\leibniz}{\Omega}        

\bmdefine{\frege}{\Lambda}         

\makeatletter
\newcommand{\tarskidsp}{\mathord%
   {\m@th\raisebox{0pt}[0pt][0pt]{$\stackrel%
   {\raisebox{-2.7pt}[0ex][0pt]{$\displaystyle \,\?\thicksim$}}%
   {\displaystyle\leibniz}$}}}
\newcommand{\tarskitxt}{\mathord%
   {\m@th\raisebox{0pt}[0pt][0pt]{$\stackrel%
   {\raisebox{-2.7pt}[0ex][0pt]{$\,\?\thicksim$}}{\displaystyle\leibniz}$}}}
\newcommand{\tarskiscr}{\mathord%
   {{\m@th\raisebox{0pt}[0pt][0pt]{$\stackrel%
   {\raisebox{-2.4pt}[0ex][0pt]{$\scriptstyle \,\?\thicksim$}}%
   {\scriptstyle\leibniz}$}}}}
\newcommand{\tarskiscrscr}{\mathord%
   {{\m@th\raisebox{0pt}[0pt][0pt]{$\stackrel%
   {\raisebox{-2pt}[0ex][0pt]{$\scriptscriptstyle \,\?\thicksim$}}%
   {\scriptscriptstyle\leibniz}$}}}}
\newcommand{\tarski}{\@ifnextchar ^ %
   {\mathchoice{\tarskidsp\kern-.07em}{\tarskitxt\kern-.07em}%
   {\tarskiscr\kern-.07em}{\tarskiscrscr\kern-.07em}}%
   {\mathchoice{\tarskidsp}{\tarskitxt}{\tarskiscr}{\tarskiscrscr}}}
\makeatother

\theoremstyle{theorem}
\newtheorem{Theorem}{Theorem}[section]
\newtheorem{Proposition}[Theorem]{Proposition}
\newtheorem{Lemma}[Theorem]{Lemma}
\newtheorem{Corollary}[Theorem]{Corollary}

\theoremstyle{definition}
\newtheorem{law}[Theorem]{Definition}
\newtheorem{exa}[Theorem]{Example}
\newtheorem{Fact}{Fact}

\theoremstyle{remark}

\newtheorem{problem}{\bf Problem}
\newtheorem{Remark}[Theorem]{Remark}

\newcommand{\C}{\boldsymbol{C}} 

\begin{document}
\title[The poset of all logics I: interpretations and lattice structures]{The poset of all logics I: \\ interpretations and lattice structure}

\author{R. Jansana and T. Moraschini}
\address{Department of Philosophy, Faculty of Philosophy, University of Barcelona, Carrer de Montalegre $6$, $08001$, Barcelona, Spain}
\email{jansana@ub.edu}
\address{Institute of Computer Science, Academy of Sciences of Czech Republic, Pod Vod\'arenskou v\v{e}\v{z}\'{i} $271/2$, $182$ $07$ Prague $8$, Czech Republic}
\email{moraschini@cs.cas.cz}
\date{\today}

\maketitle


\begin{abstract}
A notion of interpretation between arbitrary logics is introduced, and the poset $\class{Log}$ of all logics ordered under interpretability is studied. It is shown that in $\class{Log}$ infima of arbitrarily large sets exist, but binary suprema in general do not. On the other hand, the existence of suprema of sets of equivalential logics is established. The relations between $\class{Log}$ and  the lattice of interpretability types of varieties are investigated.
\end{abstract}

\section{Introduction}

Universal algebra \cite{Be11g,BuSa00} and abstract algebraic logic \cite{Cz01,AAL-AIT-f} are two disciplines that study, respectively, general algebraic structures and  propositional logics. One of their main achievements is the development of two parallel taxonomies, one of varieties (a.k.a.\ equational classes) of algebras, and the other one  of propositional logics.

More precisely, the \textit{Maltsev hierarchy} of universal algebra is a classification of varieties in terms of syntactic principles (called \textit{Maltsev conditions}) intended to describe the structure of the congruence lattices of algebras \cite{HoMcKe88,KeKi06,Pixley15can,Taylor77,Wille70}. The first, and perhaps most celebrated, example of a Maltsev condition is the requirement that a variety $\class{K}$ is congruence permutable, equivalent to the syntactic requirement of the existence of a \textit{minority term} for $\class{K}$ \cite{AIM54}, i.e.\ a ternary term $\varphi(x, y, z)$ such that
\[
\class{K} \vDash \varphi(x, x, y) \thickapprox y \thickapprox \varphi(y, x, x).
\]

Similarly, in abstract algebraic logic, the \textit{Leibniz hierarchy} is a taxonomy of propositional logics in terms of rule schemata (here called \textit{Leibniz conditions}) whose aim is to govern the interplay between lattices of deductive filters (a.k.a.\ theories) of logics and lattices of congruences of algebras \cite{BP86,BP89,Cz01,CzJa00,TMo15,Ra06a}. One of the most fundamental examples of a Leibniz condition is the requirement that a logic $\vdash$ possesses a set $\Delta(x, y)$ of binary formulas satisfying the rules
\[
\emptyset \rhd \Delta(x, x) \; \text{ and } \; x, \Delta(x, y) \rhd y,
\]
which generalize the behavior of most implication connectives. 
This requirement is equivalent to the property that the Leibniz operator of the logic $\vdash$ is monotone \cite{BP86}.

From this point of view, it is natural to wonder whether the Maltsev and Leibniz hierarchies are two faces of the same coin (see for instance \cite{JGRa11}). In a series of papers of which this one is the first (and whose next parts are \cite{JaMor19-2,JaMor19-3}) we show that this is indeed the case. More in detail, it turns out that the Maltsev hierarchy is a sort of \textit{finitary companion} of the Leibniz hierarchy of the two-deductive systems \cite{BP92}, i.e.\ substitution-invariant consequence relations between pairs of terms  understood as equations.  One of the main obstacles to establish this result is that, while there exists a precise definition of  Maltsev condition, this is not the case for what concerns Leibniz conditions (which until now were recognized on empirical grounds only). 

To clarify the notion of a Leibniz condition, we adopt an order-theoretic perspective inspired by the theory of the Maltsev hierarchy, in which varieties are ordered by means of the existence of interpretations between them \cite{Neum74,Tay73} (see also \cite{Lawvere63}). A variety $\class{K}$ is said to be \textit{interpretable} \cite{Tay73} into another variety $\class{V}$, when $\class{V}$ is term-equivalent to some variety $\class{V}^{\ast}$ whose reducts (in a smaller signature) belong to $\class{K}$. When this is true we write ${\class{K}} \leq {\class{V}}$. For instance, the variety of distributive lattices is interpretable into the one of Boolean algebras, while the variety of sets (lacking non-trivial operations) is interpretable in any variety. 
It is clear that the interpretability relation $\leq$ is a preorder on the collection of all varieties. More interestingly, the poset $\class{Var}$ associated with $\leq$ happens to be a lattice, sometimes called the lattice of \textit{interpretability types} of varieties \cite{GaTa84,Neum74}. The study of the lattice $\class{Var}$ allowed to identify the classes of models of Maltsev conditions with the filters of $\class{Var}$ that are generated by finitely presentable varieties \cite{BaBe77,Gratzer70Mal,Neum74,Tay73}.

As we mentioned, we will export this order-theoretic perspective to the realm of propositional logics that, when ordered under a suitable notion of interpretability, form the \textit{poset of all logics} $\class{Log}$. Accordingly, the aim of this paper is to describe the structure of the poset $\class{Log}$, which will be exploited to define and investigate Leibniz conditions in general in \cite{JaMor19-2,JaMor19-3}. The main results of this paper can be summarized as follows. First we establish that $\class{Log}$ is a set-complete meet-semilattice in which binary joins may fail to exist (Theorems \ref{Thm:Infima} and \ref{Thm:NotJoins}). Then we show that the proper submeet-semilattice $\class{Equiv}$ of $\class{Log}$, whose elements are equivalential logics,  happens to have joins and to be a set-complete lattice (Theorem \ref{Thm:ProtoalgebraicLattice}). We conclude by investigating the bottom and the top parts of $\class{Log}$ and by comparing the poset of all logics $\class{Log}$ with the lattice of interpretability types of varieties $\class{Var}$.

%
%

\section{Propositional logics}

For general informations on abstract algebraic logic, we refer the reader to \cite{BP86,BP89,BP92,Cz01,AAL-AIT-f,FJa09,FJaP03b,JaMoRa17,W88}. We fix a proper class of (propositional) variables $\{ x_{\alpha} \colon \alpha \in  \mathsf{OR} \}$  indexed in a one-to-one way by the ordinals. Given an algebraic language $\LL$ (from now on, simply \textit{a language}), and an infinite cardinal $\kappa$, we denote by $Fm_{\LL}(\kappa)$ the set of formulas of $\LL$ with variables among $\{x_{\alpha} \colon \alpha < \kappa \}$, and by $\Fm_{\LL}(\kappa)$ the corresponding algebra. When the language $\LL$ is clear from the context, we simply write $Fm(\kappa)$. For the sake of simplicity, we assume that  languages have no nullary operation\footnote{In the appendix we explain why this assumption is harmless, and how we can modify our approach to cover logics in languages with constants as well. However, this comes at the cost of distinguishing cases so frequently that the exposition would turn unnecessarily cumbersome.}. Note that the cardinality of $Fm_{\LL}(\kappa)$ is the maximum of  $\kappa$ and the  cardinality of $\LL$.

A \textit{logic} $\vdash$ is then a consequence relation on the set $Fm_{\LL}(\kappa)$, for some language $\LL$ and infinite cardinal $\kappa$, that is substitution invariant in the sense that for every \textit{substitution} $\sigma$ on $Fm_{\LL}(\kappa)$ and $\Gamma \cup \{ \varphi \} \subseteq Fm_{\LL}(\kappa)$,
\[
\text{if }\Gamma \vdash \varphi \text{, then }\sigma[\Gamma] \vdash \sigma(\varphi).
\]
Given a logic $\vdash$, we denote by $\LL_{\vdash}$ (resp.\ $\kappa_{\vdash}$) the language (resp.\ the cardinality of the set of variables) in which $\vdash$ is formulated. Moreover, we write $Fm(\vdash)$ as a shorthand for $Fm_{\LL_{\vdash}}(\kappa_{\vdash})$. A \textit{theorem} of $\vdash$ is  a formula $\varphi$ such that $\emptyset \vdash \varphi$.

Given an algebra $\A$ and a logic $\vdash$ in the same language, a set $F \subseteq A$ is said to be a \textit{deductive filter} of $\vdash$ on $\A$ when for every $\Gamma \cup \{ \varphi\} \subseteq Fm(\vdash)$ such that $\Gamma \vdash \varphi$ and every homomorphism $h \colon \Fm(\vdash) \to \A$, if $h[\Gamma] \subseteq F$, then $h(\varphi) \in F$. The set of deductive filters of $\vdash$ on $\A$ is a closure system, whose closure operator is denoted by $\textup{Fg}_{\vdash}^{\A}(\cdot) \colon \mathcal{P}(A) \to \mathcal{P}(A)$. Given $X \cup \{ a \} \subseteq A$, we write $\textup{Fg}_{\vdash}^{\A}(X, a)$ as a shorthand for $\textup{Fg}_{\vdash}^{\A}(X \cup \{ a\})$. Given an algebra $\B$, we also write $\B \subseteq \A$ when $\B$ is  a subalgebra of $\A$, and   $\B \leq \A$ when $\B$ is isomorphic to a subalgebra of $\A$.

\begin{Lemma}\label{Lem:compactness}
Let $\vdash$ be a logic formulated on $Fm(\kappa)$ and $\A$ an algebra.
\benroman
\item If $\B \subseteq \A$ and $X \subseteq A$, then $\textup{Fg}_{\vdash}^{\B}(X \cap B) \subseteq \textup{Fg}_{\vdash}^{\A}(X)$.
\item Let $\lambda \geq \vert Fm(\kappa) \vert$ and $X \cup Z \cup \{ a \} \subseteq A$ be such that $\vert Z \vert \leq \lambda$. If $a \in \textup{Fg}_{\vdash}^{\A}(X)$, then there is an algebra $\B \subseteq \A$ such that $\vert B \vert \leq \lambda$, $Z \subseteq B$, and $a \in \textup{Fg}_{\vdash}^{\B}(X \cap B)$.
\eroman
\end{Lemma}

\begin{proof}
Condition (i) is straightforward. Hence we detail only the proof of (ii).

 It is well-known that
\[
\textup{Fg}_{\vdash}^{\A}(X) = \bigcup_{\alpha < \lambda^{+}}V_{\alpha}
\]
where the various $V_{\alpha}$ are defined in the following way. First we set $V_{0} \coloneqq X$, and at limit ordinals we take unions. At successor ordinals we proceed as follows. If $\alpha < \lambda$, then
\begin{align*}
V_{\alpha+1} \coloneqq V_{\alpha} \cup \{ & c \in C \colon  c = f(\varphi)\text{ for some homomorphism }f \colon \Fm(\kappa) \to \A \\
&\text{and }\Gamma \cup \{ \varphi \} \subseteq Fm(\kappa) \text{ such that }\Gamma \vdash \varphi \text{ and } v[\Gamma] \subseteq V_{\alpha}\}.
\end{align*}

We claim that for every $\alpha < \lambda^{+}$ and $b \in V_{\alpha}$, there is an algebra $\B[b, \alpha] \leq \A$ such that $\vert B[b, \alpha] \vert \leq \lambda$, $Z \subseteq B[b, \alpha]$, and $b \in \textup{Fg}_{\vdash}^{\B[b, \alpha]}(X \cap B[b, \alpha])$. To prove this, we reason by induction on $\alpha \leq \lambda^{+}$. In the case where $\alpha = 0$ we take the subalgebra of $\A$ generated by $X \cup Z$. If $\alpha$ is a limit ordinal and $b \in V_{\alpha}$, then $b \in V_{\beta}$ for some $\beta < \alpha$. Therefore, with an application of the inductive hypothesis, we are done. 

Then we consider the case where $\alpha = \beta +1$. Since $b \in V_{\beta + 1}$, there are a homomorphism $f \colon \Fm(\kappa) \to \A$ and $\Gamma \cup \{ \varphi \} \subseteq Fm(\kappa)$ such that $\Gamma \vdash \varphi$, $b = f(\varphi)$, and $v[\Gamma] \subseteq V_{\beta}$. Now, for every $\gamma \in \Gamma$, we consider the algebra $\B[ \? f(\gamma), \beta]$ given by the inductive hypothesis. Let $\B[b, \alpha]$ be the subalgebra of $\A$ generated by 
\[
Z \cup f[Fm(\kappa)] \cup \bigcup_{\gamma \in \Gamma} B[\? f(\gamma), \beta ].
\]
The fact that $\vert Z \vert + \vert Fm(\kappa) \vert \leq \lambda$, and that $\vert B[\? f(\gamma), \beta] \vert \leq \lambda$ for every $\gamma \in \Gamma$ ensures that $\vert B[b, \alpha] \vert \leq \lambda$. 

It only remains to show that $b \in \textup{Fg}_{\vdash}^{\B[b, \alpha]}(X \cap B[b, \alpha])$. To this end, consider $\gamma \in \Gamma$. By the inductive hypothesis and condition (i) we obtain that
\begin{equation}\label{Eq:compactness-1}
f(\gamma) \in \textup{Fg}_{\vdash}^{\B[\? f(\gamma), \beta]}(X \cap B[\? f(\gamma), \beta]) \subseteq \textup{Fg}_{\vdash}^{\B[b, \alpha]}(X \cap B[b, \alpha]).
\end{equation}
Since $f[Fm(\kappa)] \subseteq B[b, \alpha]$, the homomorphism $f \colon \Fm(\kappa) \to \B[b, \alpha]$ is well defined. Together with the fact that $\Gamma \vdash \varphi$ and that $f[\Gamma] \subseteq \textup{Fg}_{\vdash}^{\B[b, \alpha]}(X \cap B[b, \alpha])$ by (\ref{Eq:compactness-1}), this implies that $b = f(\varphi) \in \textup{Fg}_{\vdash}^{\B[b, \alpha]}(X \cap B[b, \alpha])$, as desired. This establishes the claim.

Together with the fact that $\textup{Fg}_{\vdash}^{\A}(X) = \bigcup_{\alpha < \lambda^{+}}V_{\alpha}$, the claim concludes the proof.
\end{proof}

Given an algebra $\A$, we denote by $\Con \A$ its congruence lattice. Moreover, a congruence $\theta\in \Con \A$ is \textit{compatible} with a set $F \subseteq A$ when for every $a, b\in A$,
\[
\text{if }\langle a, b \rangle \in \theta \text{ and }a \in F\text{, then }b \in F.
\]
The \textit{Leibniz congruence} $\leibniz^{\A}F$ of $F$ is the largest congruence on $\A$ compatible with $F$. Similarly, given a logic $\vdash$ (in the same language as $\A$), we set
\[
\tarski^{\A}_{\vdash}F \coloneqq \bigcap \{ \leibniz^{\A} G \colon \text{$G$ is a deductive filter of $\vdash$ on $\A$, and }F \subseteq G \}.
\]
The relation $\tarski^{\A}_{\vdash}F$ is often called the \textit{Suszko congruence} of $F$. The congruences $\leibniz^{\A} F$ and $\tarski^{\A}_{\vdash}F$ can de described as follows \cite[Thms.\ 4.23 and 5.32]{AAL-AIT-f}:

\begin{Proposition}\label{Prop:Polynomial}
Let $\vdash$ be a logic, $\A$ an algebra, $F \subseteq A$, and $a, b\in A$.
\benroman
\item $\langle a, b \rangle \in \leibniz^{\A}F\Longleftrightarrow (p(a)\in F\textrm{ if and only if }p(b)\in F)$, for every unary polynomial function $p$ of $\A$.
\item $\langle a, b \rangle \in \tarski_{\vdash}^{\A}F\Longleftrightarrow \textup{Fg}_{\vdash}^{\A}(F, p(a))=\textup{Fg}_{\vdash}^{\A}(F, p(b))$, for every unary polynomial function $p$ of $\A$.
\eroman
\end{Proposition}

A \textit{matrix} is a pair $\langle \A, F\rangle$ such that $\A$ is an algebra and $F \subseteq A$. A matrix $\langle \A, F \rangle$ is said to be \textit{reduced} when $\leibniz^{\A}F$ is the identity relation. Moreover, we set
\[
\langle \A, F \rangle^{\ast} \coloneqq \langle \A / \leibniz^{\A}F, F/ \leibniz^{\A}F\rangle.
\]
Similarly, given a class of matrices $\class{K}$, we set
\[
\RRR(\class{K}) \coloneqq  \III \{ \langle \A, F\rangle^{\ast} \colon \langle \A, F \rangle \in \class{K} \},
\]
where $\III$ is the class operator of closing under isomorphic copies. A matrix $\langle \A, F \rangle$ is said to be \textit{trivial} when $\A$ is the trivial algebra (which we denote by ${\bf 1}$) and $F = \{ 1 \}$.

The logic \textit{induced} by a class of similar matrices $\mathsf{K}$ in $\kappa$ variables is the consequence relation $\vdash$ on $Fm(\kappa)$ defined for every $\Gamma \cup \{ \varphi \} \subseteq Fm(\kappa)$ as follows:
\begin{align*}
\Gamma \vdash \varphi \Longleftrightarrow&\text{ for every $\langle \A, F \rangle \in \class{K}$ and homomorphism }h \colon \Fm(\kappa) \to \A, \\
 &\text{ if }h[\Gamma] \subseteq F\text{, then }h(\varphi)\in F.
\end{align*}
A matrix $\langle \A, F \rangle$ is said to be a \textit{model} of a logic $\vdash$ (in the same language as $\A$) when $F$ is a deductive filter of $\vdash$ on $\A$. We set
\begin{align*}
\Mod(\vdash) &\coloneqq \{ \langle \A, F \rangle \colon \langle \A, F \rangle \text{ is a model of }\vdash\}\\
\ModS(\vdash) &\coloneqq \{ \langle \A, F \rangle \in \Mod(\vdash) \colon \tarski_{\vdash}^{\A}F \text{ is the identity relation}\}.
\end{align*}
Observe that $\vdash$ is the logic induced both by $\Mod(\vdash)$ and $\ModS(\vdash)$ \cite[Thm.\ 4.16]{AAL-AIT-f}.

We denote by  $\SSS, \PPP, \PSD$ and $\PRk$ the class operators for substructures, direct products, subdirect products, and reduced products over $\kappa$-complete filters. We assume that their application produces classes closed under isomorphic copies. Moreover, we assume that the product-style operators, when applied to empty sets of indexes, produce trivial matrices. We also consider the following class operator: given a class of matrices $\class{K}$ and an infinite cardinal $\kappa$, we define
\[
\UUU_{\kappa}(\class{K}) \coloneqq \{ \langle \A, F \rangle \colon \langle \B, F \cap B \rangle \in \class{K}\text{ for every $\kappa$-generated $\B \leq \A$}\}.
\]

\begin{Lemma}\label{Lem:PSD-Suszko}
If $\vdash$ is a logic, then $\ModS(\vdash) = \PSD\RRR(\Mod(\vdash))$.
\end{Lemma}
\begin{proof}
See \cite[Thm.\ 5.3]{Cze03}.
\end{proof}

The first equality of the following result is taken from \cite{DeJa96,Prenosil15}, and generalizes a previous result in \cite{Cz80}. 

\begin{Theorem}\label{Lem:Modstar}
Let $\class{K}$ be a class of matrices. If $\vdash$ is the logic induced by $\class{K}$ on $Fm(\kappa)$ and $\vert Fm(\kappa) \vert \leq \kappa$, then $\RRR(\Mod(\vdash)) = \RRR\SSS\PRk(\class{K})  = \RRR\UUU_{\kappa}\SSS\PPP(\class{K})$.
\end{Theorem}
\begin{proof}
Under the assumption that the cardinality of the language of a class of matrices $\class{K}$ is $\leq \kappa$, the proof of the equality $\SSS\PRk(\class{K}) = \UUU_{\kappa}\SSS\PPP(\class{K})$  is routinary.
\end{proof}

\begin{Corollary}\label{Cor:Mod-Suszko}
Let $\class{K}$ be a class of matrices. If $\vdash$ is the logic induced by $\class{K}$ on $Fm(\kappa)$ and $\vert Fm(\kappa) \vert \leq \kappa$, then $\ModS(\vdash) = \PSD\RRR\SSS\PRk(\class{K})$.
\end{Corollary}
\begin{proof}
Immediate from Lemma \ref{Lem:PSD-Suszko} and  Theorem \ref{Lem:Modstar}.
\end{proof}

\begin{Corollary}\label{Cor:Alg-general-trick}
Let $\vdash$ be the logic induced by a class of matrices $\class{K}$ on $Fm(\kappa)$. Then the algebraic reducts of the matrices in $\ModS(\vdash)$ belong to the variety generated by the algebraic reducts of the matrices in $\class{K}$.
\end{Corollary}

A logic $\vdash$ is said to be \textit{equivalential} \cite{BP86,Cz81} if there is a non-empty\footnote{In the literature the set $\Delta$ is not required to be non-empty. However, this restriction is almost immaterial as, in a fixed language, there is a unique equivalential logic with an empty $\Delta$ is the pathological \textit{almost inconsistent} logic \cite[Prop.\ 6.11.5]{AAL-AIT-f}.} set of formulas $\Delta(x, y)$ such that for every $\langle \A, F \rangle \in \Mod(\vdash)$ and $a, b \in A$,
\[
\langle a, b \rangle \in \leibniz^{\A}F \Longleftrightarrow \Delta^{\A}(a, b) \subseteq F.
\]
In this case we say that $\Delta$ is a set of \textit{congruence formulas} for $\vdash$.  Examples of equivalential logics comprise all the so-called algebraizable logics \cite{BP89}, as well as a wide range of non-algebraizable ones such as the the local consequence of the normal modal system ${\bf K}$ \cite{Ma86}.  For further information about equivalential logics, see \cite{Cz01,AAL-AIT-f,He93a,He96,He97}. 
\begin{Theorem}\label{Thm:congruence-formulas}
A logic $\vdash$ is equivalential if and only if there is a non-empty set of formulas $\Delta(x, y)$ such that for every $n$-ary connective $\ast$,
\begin{align*}
\emptyset \vdash &\Delta(x, x) \qquad x, \Delta(x, y) \vdash y\\
\bigcup_{1 \leq i \leq n}\Delta(x_{i}, y_{i}) &\vdash \Delta(\ast(x_{1}, \dots, x_{n}), \ast(y_{1}, \dots, y_{n})).
\end{align*}
In this case, $\Delta$ is a set of congruence formulas for $\vdash$ and $\ModS(\vdash) = \RRR(\Mod(\vdash))$.
\end{Theorem}

\begin{proof}
See \cite[Thms.\ 6.17 and 6.60]{AAL-AIT-f}.
\end{proof}

For equivalential logics we have the following improvement of Corollary \ref{Cor:Mod-Suszko}:

\begin{Lemma}\label{Lem:Equivalent_Alg_Sem}
Let $\vdash$ be the logic induced by a class of reduced matrices $\class{K}$ on $Fm(\kappa)$. If $\vdash$ is equivalential, then $\ModS(\vdash) = \UUU_{\kappa}\PSD\SSS(\class{K})$.
\end{Lemma}

\begin{proof}
This result is essentially \cite[Thm.\ 5.6]{RaxxNJ}.
\end{proof}

A \textit{tuple} of elements of a set $A$ is a \textit{finite} sequence of elements of $A$.

\begin{Lemma}\label{Lem:U-kappa}
If $\vdash$ is a logic on $Fm(\kappa)$ and $\lambda \geq \vert Fm(\kappa) \vert$, then $\ModS(\vdash)$ is closed under $\UUU_{\lambda}$.
\end{Lemma}
\begin{proof}
Suppose, with a view to contradiction,
 that there is a matrix $\langle \A, F \rangle \in \UUU_{\lambda}(\ModS(\vdash))$ such that $\langle \A, F \rangle \notin \ModS(\vdash)$. First observe that $\langle \A, F \rangle$ is a model of $\vdash$, since $\lambda \geq \kappa$ and $\vdash$ is defined on $Fm(\kappa)$. Then the congruence $\tarski_{\vdash}^{\A}F$ is not the identity relation. Together with Proposition \ref{Prop:Polynomial}(ii), this implies that there are two different $a, b \in A$ such that for every $\varphi(x, \vec{y}) \in Fm(\omega)$, and every tuple $\vec{c} \in A$,
\begin{equation}\label{Eq:first-equation}
\textup{Fg}_{\vdash}^{\A}(F, \varphi(a, \vec{c})) = \textup{Fg}_{\vdash}(F, \varphi^{\A}(b, \vec{c})).
\end{equation}

We define a chain (under the inclusion relation) $\langle \B_{\alpha} \colon \alpha < \lambda\rangle$ of subalgebras of $\A$ as follows. First we let $\B_{0}$ be the subalgebra of $\A$ generated by $\{ a, b \}$. At limit ordinals we take unions. Now, suppose that $\B_{\alpha}$ has already been defined and that $\alpha < \lambda$. Consider a formula $\varphi(x, \vec{y}) \in Fm(\omega)$ and a tuple $\vec{c} \in B_{\alpha}$. By Lemma \ref{Lem:compactness} and (\ref{Eq:first-equation}) there is a subalgebra $\B[\varphi, \vec{c}, \alpha] \leq \A$ such that $\vert B[\varphi, \vec{c}, \alpha] \vert \leq \lambda$, $a, b, \vec{c} \in B[\varphi, \vec{c}]$ and 
\begin{equation}\label{Eq:compactness-2}
\textup{Fg}_{\vdash}^{\B[\varphi, \vec{c}, \alpha]}(F \cap B[\varphi, \vec{c}, \alpha], \varphi(a, \vec{c})) = \textup{Fg}_{\vdash}^{\B[\varphi, \vec{c}, \alpha]}(F \cap B[\varphi, \vec{c}, \alpha], \varphi(b, \vec{c})).
\end{equation}
Then we let $\B_{\alpha+1}^{\ast}$ be the subalgebra of $\A$ generated by the union of the various $B[\varphi, \vec{c}, \alpha]$, and  $\B_{\alpha}$  the subalgebra of $\A$ generated by $B_{\alpha} \cup B_{\alpha+1}^{\ast}$. 

Now we set
\[
\langle \B, G \rangle \coloneqq \langle \bigcup_{\alpha < \lambda} \B_{\alpha}, F \cap \bigcup_{\alpha < \lambda}B_{\alpha} \rangle.
\]
Bearing in mind that $\vert B[\varphi, \vec{c}, \alpha] \vert + \vert Fm(\kappa) \vert \leq \lambda$, an easy induction shows that $\vert B_{\alpha} \vert \leq \lambda$ for every $\alpha < \lambda$. As a consequence, we obtain that $\vert B \vert \leq \lambda$ and, therefore, that $\B$ is $\lambda$-generated. Together with $\langle \A, F \rangle \in \UUU_{\lambda}(\ModS(\vdash))$, this implies that $\langle \B, G \rangle \in \ModS(\vdash)$.

Now, the fact that $\langle \B, G \rangle \in \ModS(\vdash)$ implies that $\tarski_{\vdash}^{\B}G$ is the identity relation and, therefore, that $\langle a, b \rangle \notin \tarski_{\vdash}^{\B}G$. By Lemma \ref{Prop:Polynomial}(ii) we can assume without loss of generality that there are a formula $\varphi(x, \vec{y}) \in Fm(\omega)$ and a tuple $\vec{c} \in B$ such that $\varphi(a, \vec{c}) \notin \textup{Fg}_{\vdash}^{\B}(G, \varphi(b, \vec{c}))$. Observe that there is $\alpha < \lambda$ such that $\vec{c} \in B_{\alpha}$. By (\ref{Eq:compactness-2}) and Lemma \ref{Lem:compactness}(i) we obtain that
\begin{align*}
\varphi(a, \vec{c}) &\in \textup{Fg}_{\vdash}^{\B[\varphi, \vec{c}, \alpha]}(F \cap B[\varphi, \vec{c}, \alpha], \varphi(b, \vec{c}))\\
&= \textup{Fg}_{\vdash}^{\B[\varphi, \vec{c}, \alpha]}(G \cap B[\varphi, \vec{c}, \alpha], \varphi(b, \vec{c}))\\
& \subseteq \textup{Fg}_{\vdash}^{\B}(G, \varphi(b, \vec{c})).
\end{align*}
But this contradicts the fact that $\varphi(a, \vec{c}) \notin \textup{Fg}_{\vdash}^{\B}(G, \varphi(b, \vec{c}))$. Hence we reached a contradiction, as desired.
\end{proof}

\section{Interpretations}

\begin{law}
Let $\LL$ and $\LL'$ be two languages. A \textit{translation} $\btau$ of $\LL$ into $\LL'$ is a map that associates an $n$-ary formula $\btau(\ast)$ of $\LL'$ in variables $x_{1}, \dots, x_{n}$ to every $n$-ary function symbol $\ast$ of $\LL$.
\end{law}

Let $\btau$ be a translation of $\LL$ into $\LL'$. Given two infinite cardinals $\kappa \leq \lambda$ and a formula $\varphi \in Fm_{\LL}(\kappa)$, we define a formula $\btau(\varphi) \in Fm_{\LL'}(\lambda)$ by recursion as follows. If $\varphi = x_{\alpha}$ for some $\alpha < \kappa$, then $\btau(\varphi) \coloneqq x_{\alpha}$. Moreover, if $\varphi = \ast(\psi_{1}, \dots, \psi_{n})$ for some $n$-ary function symbol $\ast$ of $\LL$, then $\btau(\varphi) \coloneqq \btau(\ast)(\btau(\psi_{1}), \dots, \btau(\psi_{n}))$. We extend this notation to sets of formulas $\Gamma \subseteq Fm_{\LL}(\kappa)$, by setting $\btau[\Gamma] \coloneqq \{ \btau(\gamma) \colon \gamma \in \Gamma \}$. Note that the variables of $\btau(\varphi)$ are among the variables in $\varphi$.

Moreover, given an $\LL'$-algebra $\A$, we let $\A^{\btau}$ be the  $\LL$-algebra, whose universe is $A$, and whose $n$-ary operations $\ast$ are interpreted as follows:
\[
\ast^{\A^{\btau}}(a_1, \ldots, a_n) \coloneqq \btau(\ast)^{\A}(a_1, \ldots, a_n), \text{ for every }a_{1}, \dots, a_{n} \in A.
\]
By induction on the construction of the formulas we obtain that for every $\varphi(z_1, \ldots, z_n)  \in  Fm_{\LL}(\kappa)$ and every $a_1, \ldots, a_n \in A$, 
\[
\btau(\varphi)^{\A}(a_1, \ldots, a_n) = \varphi^{\A^{\btau}}(a_1, \ldots, a_n).
\]

\begin{law}
\label{def:interpretation}
Let $\vdash$ and $\vdash'$ be two logics. An \textit{interpretation} of $\vdash$ into $\vdash'$ is a translation $\btau$ of $\LL_{\vdash}$ into  $\LL_{\vdash'}$ such that
\[
\text{if }\langle \A, F \rangle \in \ModS(\vdash')\text{, then }\langle \A^{\btau}, F \rangle \in \ModS(\vdash).
\]
\end{law}

For instance, for every given logic the identity map is an interpretation of it into any of its extensions.

\begin{Proposition}\label{Prop:interpretation}
If $\btau$ is an interpretation of $\vdash$ into $\vdash'$, and $\langle \A, F \rangle \in \Mod(\vdash')$, then $\langle \A^{\btau}, F \rangle \in \Mod(\vdash)$. Moreover, if $\lambda \leq \kappa_{\vdash} \leq \kappa_{\vdash'}$, then for every $\Gamma \cup \{ \varphi \} \subseteq Fm_{\LL_{\vdash}}(\lambda)$,
\[
\text{if }\Gamma \vdash\varphi\text{, then }\btau[\Gamma] \vdash' \btau(\varphi).
\]
\end{Proposition}

\begin{Proposition}\label{Prop:interpretation2}
Let $\vdash$ and $\vdash'$ be two logics and $\btau$ be a translation of $\LL_{\vdash}$ into $\LL_{\vdash'}$. Then $\btau$ is an interpretation of $\vdash$ into $\vdash'$ if and only if $\langle \A^{\btau}, F \rangle \in \ModS(\vdash)$ for every $\langle \A, F \rangle \in \RRR(\Mod(\vdash'))$.
\end{Proposition}

\begin{proof}
The ``only if'' part is immediate. The ``if'' one is a consequence of Lemma \ref{Lem:PSD-Suszko}. 
\end{proof}

When there is an interpretation of $\vdash$ into $\vdash'$ we write ${\vdash} \leq {\vdash'}$ and say that $\vdash$ is \textit{interpretable} into $\vdash'$. Similarly, we say that $\vdash$ and $\vdash'$ are \textit{equi-interpretable} if ${\vdash} \leq {\vdash'}$ and ${\vdash'} \leq {\vdash}$. Given a logic $\vdash$, we denote by $\llbracket  \vdash \rrbracket$ the class of all logics which are equi-interpretable with $\vdash$. It is clear that relation $\leq$ is a preorder on the \textit{proper class} of all logics, and that it induces a partial order on the \textit{collection} of all classes of the form $\llbracket  \vdash \rrbracket$.  The latter poset constitutes the object of study of this work.

\begin{law}
We denote by $\class{Log}$ the \textit{poset of all logics}, i.e.\ the poset whose universe is $\{ \llbracket  \vdash \rrbracket \colon {\vdash} \text{ is a logic}\}$ equipped with the partial order $\leq$, defined as follows:
\[
\llbracket  \vdash \rrbracket \leq \llbracket  \vdash' \rrbracket \Longleftrightarrow {\vdash} \leq {\vdash'}.
\]
\end{law}

\begin{Remark}
The reader may feel reassured by learning that, despite our reference to classes and collections, the results of this work can be formulated entirely in ZFC. This is because our statements can be phrased equivalently as speaking about logics ordered under the \textit{preorder} $\leq$ by modifying the statements about posets to statements about preorders in the natural way. It is therefore only for the sake of simplicity that we found convenient to work with the poset $\class{Log}$ whose elements are, strictly speaking, proper classes.
\qed
\end{Remark}

The notion of interpretability can be broken into two halves as follows:

\begin{law} Let $\vdash$ and $\vdash'$ be logics.
\benroman
\item $\vdash$ and $\vdash'$ are \textit{term-equivalent} if there are interpretations $\btau$ of $\vdash$ into $\vdash'$ and $\brho$ of $\vdash'$ into $\vdash$ such that
\[
\langle \A, F \rangle =\langle \A^{\btau\brho}, F \rangle \; \text{ and } \; \langle \B, G \rangle = \langle \B^{\brho\btau}, G \rangle
\]
for every $\langle \A, F \rangle \in \ModS(\vdash')$ and $\langle \B, G \rangle \in \ModS(\vdash)$.
\item $\vdash'$ is a \textit{compatible expansion} of $\vdash$ if $\LL_{\vdash} \subseteq \LL_{\vdash'}$ and the $\LL_{\vdash}$-reducts of the structures in $\ModS(\vdash')$ belong to $\ModS(\vdash)$.
\eroman
 \end{law}

\begin{Proposition}\label{Prop:CompatibleExpansions}
Let $\vdash$ and $\vdash'$ be logics. Then ${\vdash} \leq {\vdash'}$ if and only if $\vdash'$ is term-equivalent to a compatible expansion of $\vdash$.
\end{Proposition}

\begin{proof}
The ``if'' part is immediate. To prove the ``only if'' part, suppose that there is an interpretation $\btau$ of $\vdash$ into $\vdash'$. We can assume without loss of generality that the sets of function symbols of $\vdash$ and $\vdash'$ are disjoint. Then let $\LL$ be the language extending $\LL_{\vdash'}$  with the symbols of $\vdash$. Given a matrix $\langle \A, F \rangle \in \ModS(\vdash')$, we denote by $\A_{\LL}$ the $\LL$-algebra obtained by enriching $\A$ with the following interpretation of $n$-ary symbols $\ast$ of $\vdash$: for every $a_{1}, \dots, a_{n} \in A$,
\[
\ast^{\A_{\LL}}(a_{1}, \dots, a_{n}) \coloneqq \btau(\ast)^{\A}(a_{1}, \dots, a_{n}).
\]

Then consider the class of matrices $\class{K} \coloneqq \{ \langle \A_{\LL}, F \rangle \colon \langle \A, F \rangle \in \ModS(\vdash') \}$, and let $\vdash''$ be the logic on $Fm_{\LL}(\kappa_{\vdash'})$ induced by $\class{K}$. It is not hard to see that $\ModS(\vdash'') = \class{K}$. Together with the fact that $\btau$ is an interpretation of $\vdash$ into $\vdash'$, this implies that $\vdash''$ is a compatible expansion of $\vdash$. As it is clear that $\vdash'$ and $\vdash''$ are term-equivalent, we are done.
\end{proof}

The following  is instrumental to construct concrete interpretations.

\begin{Proposition}\label{Prop:alg-translations}
Let $\class{K}$ be a class of reduced matrices that induces an equivalential logic $\vdash'$. Moreover, let $\vdash$ be a logic such that $\kappa_{\vdash'} \geq \vert Fm(\vdash) \vert$. A translation $\btau$ of $\LL_{\vdash}$ into $\LL_{\vdash'}$ is an interpretation of $\vdash$ into $\vdash'$ if and only if $\langle \A^{\btau}, F \rangle \in \ModS(\vdash)$ for every $\langle \A, F \rangle \in \SSS(\class{K})$.
\end{Proposition}

\begin{proof}
The ``if'' part follows from the fact that $\SSS(\class{K}) \subseteq \ModS(\vdash')$ by Lemma \ref{Lem:Equivalent_Alg_Sem}.

To prove the ``only if'' part, suppose that $\langle \A^{\btau}, F \rangle \in \ModS(\vdash)$ for every $\langle \A, F \rangle \in \SSS(\class{K})$. By Lemmas \ref{Lem:PSD-Suszko} and \ref{Lem:U-kappa} this yields that $\langle \A^{\btau}, F \rangle \in \ModS(\vdash)$ for every $\langle \A, F \rangle \in \UUU_{\kappa_{\vdash'}}\PSD\SSS(\class{K})$. With an application of Lemma \ref{Lem:Equivalent_Alg_Sem}, we conclude that $\langle \A^{\btau}, F \rangle \in \ModS(\vdash)$ for every $\langle \A, F \rangle \in \ModS(\vdash')$ and, therefore, that $\btau$ is an interpretation of $\vdash$ into $\vdash'$.
\end{proof}

\section{Existence of infima of sets}

A basic question about the poset $\class{Log}$ is whether it is a lattice or not. It turns out that $\class{Log}$ has infima of arbitrarily large sets, but unfortunately may lack even finite suprema. In this section we describe a construction that supplies an explicit description of infima.


\begin{law}
Given a family $\{\LL_i \colon i \in I \}$ of languages, we denote by $\bigotimes_{i \in I}\LL_i$ the language whose $n$-ary symbols $\ast$ are sequences of the form
\[
\ast =  \langle \varphi_{i}(x_1, \ldots, x_n) \colon i \in I \rangle,
\] 
where $\varphi_{i}(x_{1}, \dots, x_{n}) \in Fm_{\LL_{i}}(\omega)$ for every $i \in I$. Keeping this in mind, consider a family $J = \{ \A_{i} \colon i \in I \}$ in which $\A_{i}$ is an $\LL_{i}$-algebra, for every $i \in I$. The \textit{non-indexed product} $\bigotimes_{i \in I}\A_i$ of $J$ is the $\bigotimes_{i \in I}\LL_i$-algebra defined as follows:
\benroman
\item the universe of $\bigotimes_{i \in I}\A_i$ is the Cartesian product $\prod_{i \in I}A_{i}$, and
\item the $n$-ary symbols $\ast =  \langle \varphi_{i}(x_1, \ldots, x_n) \colon i \in I \rangle$ are interpreted as
\[
\ast^{\bigotimes_{i \in I} \A_{i}}(\vec{a}_{1}, \dots, \vec{a}_{n}) \coloneqq \langle \varphi_{i}^{\A_{i}}(\vec{a}_{1}(i), \dots, \vec{a}_{n}(i)) \colon i \in I \rangle,
\]
\eroman
for every $\vec{a}_{1}, \dots, \vec{a}_{n} \in \prod_{i \in I}A_{i}$. \qed
\end{law}
Non-indexed products of algebras found various applications in universal algebra, especially in the theory of Maltsev conditions \cite{BaBe77,GaTa84,Gratzer70Mal,Neum74,Tay73}. We use the terminology of these papers and extend it to families of matrices and logics.

\begin{law}
The \textit{non-indexed product} of a family $\{ \langle \A_{i}, F_{i}\rangle \colon  i \in I\}$ of matrices is defined in a similar fashion, by setting
\[
\pushQED{\qed} \bigotimes_{i \in I}\langle \A_{i}, F_{i}\rangle \coloneqq \langle \bigotimes_{i \in I}\A_{i}, \prod_{i \in I}F_{i}\rangle. \qedhere \popQED
\]
\end{law}

\begin{Remark}
If  $\{ {\vdash_{i}} \colon i \in I \}$ is a family of   logics, then the cardinal of  $\bigotimes_{i \in I}{\LL_{\vdash_i}}$ is lesser than or equal to $\prod_{i \in I}\vert Fm(\vdash_{i})\vert$. Moreover, if $\kappa \geq \prod_{i \in I}\vert Fm(\vdash_{i})\vert$ and  $Fm(\kappa)$ is the set of formulas of $\bigotimes_{i \in I}{\LL_{\vdash_i}}$ in $\kappa$ variables, then 
$\vert Fm(\kappa)\vert \leq \kappa$.
\qed
\end{Remark}
 
 Given a collection $\{ \class{K}_{i} \colon i \in I \}$ in which $\class{K}_{i}$ is a class of $\LL_{i}$-matrices and $I$ is a set, we define
\[
\bigotimes_{i \in I} \class{K}_{i} \coloneqq \III \{ \bigotimes_{i \in I}\langle \A_{i}, F_{i}\rangle \colon \langle \A_{i}, F_{i}\rangle \in \class{K}_{i} \}.
\]

A submatrix $\langle \A, F \rangle \subseteq \bigotimes_{i \in I}\langle \A_{i}, F_{i}\rangle$ is said to be a \textit{non-indexed subdirect product} of $\{ \langle \A_{i}, F_{i}\rangle \colon  i \in I\}$, in symbols $\langle \A, F \rangle \subseteq_{\;\textup{sd}} \bigotimes_{i \in I}\langle \A_{i}, F_{i}\rangle$, if the projection maps $\pi_{i} \colon A \to A_{i}$ are surjective. We write $\langle \A, F \rangle \leq_{\; \textup{sd}} \bigotimes_{i \in I}\langle \A_{i}, F_{i}\rangle$ to indicate that $\langle \A, F \rangle$ is isomorphic to a matrix $\langle \B, G\rangle$ such that $\langle \B, G \rangle \subseteq_{\;\textup{sd}} \bigotimes_{i \in I}\langle \A_{i}, F_{i}\rangle$.

\begin{law}
Let $\{ \vdash_{i} \colon i \in I \}$ be a family of logics. The \textit{non-indexed product} $\bigotimes_{i \in I}{\vdash_{i}}$ of $\{ \vdash_{i} \colon i \in I \}$ is the logic in the language $\bigotimes_{i \in I}\LL_i$ formulated in $\kappa$ variables and induced by the class of matrices $\bigotimes_{i \in I}\ModS(\vdash_{i})$, where
\[
\kappa \coloneqq   \prod_{i \in I}\vert Fm(\vdash_{i}) \vert.
\]

When $I = \emptyset$, we stipulate that $\bigotimes_{i \in I} {\vdash_{i}}$ is the logic in the empty language formulated in countably many variables and induced by the trivial matrix $\langle \boldsymbol{1}, \{ 1 \}\rangle$. \qed
\end{law}

Our aim is to prove that $\llbracket\bigotimes_{i \in I}{\vdash_{i}} \rrbracket$ is the infimum of $\{ \llbracket \vdash_{i} \rrbracket \colon i \in I \}$ in $\class{Log}$. To this end, we rely on the following characterization of $\ModS(\bigotimes_{i \in I}{\vdash_{i}})$, to be established later on.

\begin{Proposition}\label{Prop:Taylor-Suszko}
If $\{ \vdash_{i} \colon i \in I \}$ is a family of logics, then 
\[
\ModS(\bigotimes_{i \in I}{\vdash_{i}}) =\PSD(\bigotimes_{i \in I} \RRR(\Mod(\vdash_{i}))) =  \PSD(\bigotimes_{i \in I} \ModS(\vdash_{i})).
\]
Moreover,
\[
\ModS(\bigotimes_{i \in I}{\vdash_{i}}) = \{\langle \A, F\rangle:\langle \A, F\rangle \leq_{\textup{sd}}  \bigotimes_{i \in I}\langle \A_{i}, F_{i}\rangle \text{ for some  } \langle \A_{i}, F_{i}\rangle\in  \ModS(\vdash_i)\}.
\]

\end{Proposition}

As we promised, we obtain the following:

\begin{Theorem}\label{Thm:Infima}
The infimum of a set $\{ \llbracket \vdash_{i} \rrbracket \colon i \in I \} \subseteq \class{Log}$ is $\llbracket \bigotimes_{i \in I}{\vdash_{i}} \rrbracket$. Thus $\class{Log}$ is a set-complete meet-semilattice, i.e.\ infima of subsets of $\class{Log}$ exist. 

\end{Theorem}

\begin{proof}
First we show that ${\bigotimes_{i \in I}{\vdash_{i}}}\leq {\vdash_{j}}$ for every $j \in I$. To this end, consider the map $\btau$ that sends every $n$-ary basic operation of $\bigotimes_{i \in I}{\vdash_{i}}$ to its $j$-th component (which is an $n$-ary  term of $\vdash_{j}$). Consider $\langle \A, F \rangle \in \ModS(\vdash_{j})$. It is clear that $\langle \A^{\btau}, F \rangle \cong \bigotimes_{i \in I}\langle \A_{i}, F_{i}\rangle$ such that $\langle \A_{i}, F_{i}\rangle$ is the trivial $\mathscr{L}_{i}$-matrix for every $i \in I \smallsetminus \{ j \}$, and $\langle \A_{j}, F_{j} \rangle \coloneqq \langle \A, F \rangle$. By Proposition \ref{Prop:Taylor-Suszko} we have 
\[
\langle \A^{\btau}, F \rangle \, \cong\,  \bigotimes_{i \in I}\langle \A_{i}, F_{i}\rangle \in \bigotimes_{i \in I}\ModS(\vdash_{i}) \subseteq \ModS(\bigotimes_{i \in I}{\vdash_{i}}).
\]
In particular, this means that $\btau$ is an interpretation of $\bigotimes_{i \in I}{\vdash_{i}}$ in $\vdash_{j}$, thus
$\bigotimes_{i \in I}{\vdash_{i}} \leq {\vdash_{j}}$. 
As a consequence, $\llbracket \bigotimes_{i \in I}{\vdash_{i}} \rrbracket$ is a lower bound of $\{ \llbracket \vdash_{i} \rrbracket \colon i \in I \}$.

To prove that $\{ \llbracket \vdash_{i} \rrbracket \colon i \in I \}$ is the greatest lower bound of $\llbracket \bigotimes_{i \in I}{\vdash_{i}} \rrbracket$, consider a logic $\vdash$ such that ${\vdash} \leq {\vdash_{i}}$ for every $i \in I$. Then for each $i \in I$ there is an interpretation $\btau_{i}$ of $\vdash$ into $\vdash_{i}$. Let $\btau$ be the map that associates with every basic $n$-ary symbol $\ast$ of $\vdash$ the following $n$-ary term of $\bigotimes_{i \in I}{\vdash_{i}}$:
\[
\btau(\ast) \coloneqq \langle \btau_{i}(\ast) \colon i \in I \rangle.
\]
Now, consider a matrix $\langle \A, F \rangle \in \ModS(\bigotimes_{i \in I}{\vdash_{i}})$. From Proposition \ref{Prop:Taylor-Suszko} it follows that $\langle \A, F \rangle \leq \prod_{j \in J} (\bigotimes_{i \in I}\langle \A_{i}^{j}, F_{i}^{j}\rangle)$ is a subdirect product for some $\langle \A_{i}^{j}, F_{i}^{j}\rangle \in \ModS(\vdash_{i})$. It is easy to see that
\[
\langle \A^{\btau}, F \rangle \; \leq\;  \prod_{j \in J} \Big(\prod_{i \in I}\langle (\A_{i}^{j})^{\btau_{i}}, F_{i}^{j}\rangle \Big)
\]
is also a subdirect product. Since each $\btau_{i}$ is an interpretation of $\vdash$ into $\vdash_{i}$, we conclude that 
\[
\langle \A^{\btau}, F \rangle \in \PSD\PPP(\ModS(\vdash)) = \PSD(\ModS(\vdash)).
\]
Together with the fact that $\ModS(\vdash)$ is closed under subdirect products  by Lemma \ref{Lem:PSD-Suszko}, this yields that $\langle \A^{\btau}, F \rangle \in \ModS(\vdash)$. Hence we conclude that $\bigotimes_{i \in I}{\vdash_{i}} \leq {\vdash}$.
\end{proof}



The remaining part of this section is devote to prove Proposition \ref{Prop:Taylor-Suszko}. The proof proceeds through a series of technical observations. 

%

\begin{Lemma}\label{Lem:Taylor-Leibniz-a}
If $\langle \A, F \rangle \subseteq_{\textup{sd}} \bigotimes_{i \in I}\langle \A_{i}, F_{i}\rangle$ and $F \ne \emptyset$, then for every $\vec{a}, \vec{c} \in A$,
\[
\langle \vec{a}, \vec{c} \? \rangle \in \leibniz^{\A}F  \Longleftrightarrow \text{ for every }i \in I, \langle \vec{a}(i), \vec{c}(i)\rangle \in \leibniz^{\A_{i}}F_{i}.
\]
\end{Lemma}
\begin{proof}
The right-to-left direction is an easy exercise.\ To prove the left-to-right direction, suppose that $\langle \vec{a}, \vec{c}\rangle \in \leibniz^{\A}F$. By Lemma \ref{Prop:Polynomial}(i), given an arbitrary $j \in I$, we need to show that $p(\vec{a}(j)) \in F$ iff $p(\vec{c}(j)) \in F$, for every unary polynomial function $p(x)$ of $\A_{j}$. To this end, consider a formula  $\varphi(x, y_{1}, \dots, y_{n})$ of $\A_{j}$ and elements $e_{1}, \dots, e_{n} \in A_{j}$ such that
\begin{equation}\label{Eq:Taylor-trick1}
\varphi^{\A_{j}}(\vec{a}(j), e_{1}, \dots, e_{n}) \in F_{j}.
\end{equation}
Since $\pi_{j} \colon A \to A_{j}$ is surjective, there are $\vec{e}_{1}, \dots, \vec{e}_{n} \in A$ whose $j$-th components are respectively $e_{1}, \dots, e_{n}$. Moreover, as $F \ne \emptyset$, we can choose an element $\vec{e} \in F$. Then consider the basic operation 
\[
\psi(x, y_{1}, \dots, y_{n}, z) \coloneqq \langle \psi_{i}(x, y_{1}, \dots, y_{n}, z) \colon i \in I \rangle
\]
of $\A$, where $\psi_{j} = \varphi$, and $\psi_{i} = z$ for every $i \in J \smallsetminus \{ j \}$. We have that for every $i \in I$,
\[
\psi(\vec{a}, \vec{e}_{1}, \dots, \vec{e}_{n}, \vec{e})(i)= \left\{ \begin{array}{ll}
\varphi^{\A_{j}}(\vec{a}(j), e_{1}, \dots, e_{n}) & \text{if $i= j$}\\
  \vec{e}(i) & \text{otherwise.}\\
  \end{array} \right.
\]
Together with (\ref{Eq:Taylor-trick1}) and $\vec{e} \in F$, this implies that $\psi(\vec{a}, \vec{e}_{1}, \dots, \vec{e}_{n}, \vec{e}) \in F$. Since $\langle \vec{a}, \vec{c }\?\rangle \in \leibniz^{\A}F$, we obtain that $\psi(\vec{c}, \vec{e}_{1}, \dots, \vec{e}_{n}, \vec{e}) \in F$ as well. In particular, this means that
\[
\varphi^{\A_{j}}(\vec{c}(j), e_{1}, \dots, e_{n}) = \psi(\vec{c}, \vec{e}_{1}, \dots, \vec{e}_{n}, \vec{e})(j) \in F_{j}.
\]
Hence we conclude that $\langle \vec{a}(j), \vec{c}(j)\rangle \in \leibniz^{\A_{j}}F_{j}$,  as desired.
\end{proof}

\begin{Corollary}\label{Cor:Taylor-Leibniz}
If $\langle \A, F \rangle \subseteq_{\textup{sd}} \bigotimes_{i \in I}\langle \A_{i}, F_{i}\rangle$ and $F \ne \emptyset$, then
\benroman
\item if the matrices in $\{\langle \A_{i}, F_{i}\rangle: i \in I\}$ are reduced, then so is $\langle \A, F \rangle$;
\item $\langle \A, F \rangle^{\ast}  \leq_{\textup{sd}}\bigotimes_{i \in I}\langle\A_{i}, F_{i}\rangle^{\ast}$.
\eroman
\end{Corollary}
\begin{proof}
Condition (i) is an immediate consequence of Lemma \ref{Lem:Taylor-Leibniz-a}. 
To prove condition (ii), consider the map $f \colon \langle \A, F \rangle^{\ast} \to \bigotimes_{i \in I}\langle\A_{i}, F_{i}\rangle^{\ast}$ defined as
\[
f(a/\leibniz^{\A}F) \coloneqq \langle a(i)/\leibniz^{\A_{i}}F_{i} \colon i \in I\rangle
\]
for every $a \in A$. From Lemma \ref{Lem:Taylor-Leibniz-a} it follows that $f$ is a well-defined embedding. Together with the fact that $\langle \A, F \rangle \subseteq_{\textup{sd}} \bigotimes_{i \in I}\langle \A_{i}, F_{i}\rangle$, this implies that $\langle \A, F \rangle^{\ast}  \leq_{\textup{sd}}\bigotimes_{i \in I}\langle\A_{i}, F_{i}\rangle^{\ast}$.
\end{proof}

\begin{Proposition}\label{Prop:Taylor-thms}
Let $\{ \vdash_{i} \colon i \in I \}$ be a family of logics. The logic $\bigotimes_{i \in I}{\vdash_{i}}$ has theorems if and only if each $\vdash_{i}$ has theorems. 
\end{Proposition}

\begin{proof}
The ``only if'' part is immediate. To prove the ``if'' one, suppose that each $\vdash_{i}$ has a theorem $\varphi_{i}$. By substitution invariance, we can assume that $\varphi_{i} = \varphi_{i}(x)$. Then the formula $\varphi(x) \coloneqq \langle \varphi_{i}(x) \colon i \in I \rangle$ is a theorem of $\bigotimes_{i \in I}{\vdash_{i}}$.
\end{proof}

\begin{Lemma}\label{Lem:Taylor-reduced}
Let $\{ \vdash_{i} \colon i \in I \}$ be a family of logics, and $\langle \A, F \rangle$ a matrix such that $F \ne \emptyset$. The following conditions are equivalent:
\benroman
\item $\langle \A, F \rangle \in \RRR(\Mod(\bigotimes_{i \in I}{\vdash_{i}}))$.
\item $\langle \A, F \rangle \leq_{\textup{sd}} \bigotimes_{i \in I}\langle \A_{i}, F_{i} \rangle$, for some $\langle \A_{i}, F_{i}\rangle \in \RRR(\Mod(\vdash_{i}))$.
\eroman
\end{Lemma}

\begin{proof}
(i)$\Rightarrow$(ii): Let $\kappa \coloneqq   \prod_{i \in I}\vert Fm(\vdash_{i}) \vert$ and $Fm(\kappa)$  the set of formulas of $\bigotimes_{i \in I} {\vdash_{i}}$ in $\kappa$ variables. We know that $\kappa \geq \vert Fm(\kappa) \vert$. Since  $\bigotimes_{i \in I} \vdash_{i}$ is the logic on $Fm(\kappa)$ induced by $\bigotimes_{i \in I}\ModS(\vdash_{i})$, we can apply Theorem \ref{Lem:Modstar} yielding
\[
\langle \A, F \rangle \in \RRR\SSS\PRk(\bigotimes_{i \in I}\ModS(\vdash_{i})).
\]

Then there are a matrix $\langle \B, G \rangle$, a family of matrices $\{ \langle \B_{i}^{j}, G_{i}^{j}\rangle \colon i \in I, j \in J \}$, and a $\kappa^{+}$-complete filter $F$ on $J$ such that $\langle \B, G \rangle^{\ast} = \langle \A, F \rangle$, $\langle\B_{i}^{j}, G_{i}^{j}\rangle \in \ModS(\vdash_{i})$, and
\begin{equation}\label{Eq:the-beta-map}
\langle \B, G \rangle \leq \Big(\prod_{j \in J} (\bigotimes_{i \in I} \langle \B_{i}^{j}, G_{i}^{j}\rangle) \Big)/ F.
\end{equation}
It is easy to see that the map
\[
f \colon \prod_{j \in J} (\bigotimes_{i \in I} \langle \B_{i}^{j}, G_{i}^{j}\rangle) \to  \bigotimes_{i \in I} (\prod_{j \in J}\langle \B_{i}^{j}, G_{i}^{j}\rangle),
\]
defined by the rule
\[
f(\vec{a})(i)(j) \coloneqq \vec{a}(j)(i), \text{ for every }i \in I, j \in J,
\]
is an isomorphism. We shall see that also the map
\[
g \colon \prod_{j \in J} (\bigotimes_{i \in I} \langle \B_{i}^{j}, G_{i}^{j}\rangle) / F \to \bigotimes_{i \in I} (\prod_{j \in J}\langle \B_{i}^{j}, G_{i}^{j}\rangle / F),
\]
defined by the rule 
\[
g(\vec{a}/ F)(i) \coloneqq f(\vec{a})(i)/F \text{, for every }i \in I,
\]
is an isomorphism.\ The proof that $g$ is a well-defined surjective homomorphism is routinary. To prove that $g$ is also injective, consider  $\vec{a}, \vec{c} \in \prod_{j \in J} (\bigotimes_{i \in I} \langle \B_{i}^{j}, G_{i}^{j}\rangle) $ such that $g(\vec{a}/ F) = g(\vec{c}/ F)$, i.e.\ that $f(\vec{a})(i) / F = f(\vec{c}) / F$ for every $i \in I$. Since $\kappa \geq \vert I \vert$ and $F$ is $\kappa^{+}$-complete, we have 
\begin{align*}
\{ j \in J \colon \vec{a}(j) = \vec{c}(j) \} &= \bigcap_{i \in I} \{ j \in J \colon \vec{a}(j)(i) = \vec{c}(j)(i) \}\\
& = \bigcap_{i \in I} \{ j \in J \colon f(\vec{a})(i)(j) = f(\vec{c})(i)(j) \}\\
& \in F.
\end{align*}
Hence $\vec{a}/ F = \vec{c}/ F$ and, therefore, $g$ is injective. This establishes that $g$ is an isomorphism.

Together with (\ref{Eq:the-beta-map}), this yields that $\langle \B, G \rangle \leq \bigotimes_{i \in I} (\prod_{j \in J}\langle \B_{i}^{j}, G_{i}^{j}\rangle / F)$. As a consequence, there are $\langle \A_{i}, F_{i}\rangle \in \SSS\PRk(\ModS(\vdash_{i}))$ such that $\langle \B, G \rangle \leq_{\textup{sd}} \bigotimes_{i \in I}\langle \A_{i}, F_{i}\rangle$. Together with Corollary \ref{Cor:Taylor-Leibniz}, this implies that
\[
\langle \A, F\rangle \leq_{\; \textup{sd}} \bigotimes_{i \in I}\langle \A_{i}, F_{i}\rangle^{\ast}
\]
where $\langle \A_{i}, F_{i} \rangle^{\ast} \in \RRR\SSS\PRk(\ModS(\vdash_{i}))$. As $\kappa \geq \vert Fm(\vdash_{i})\vert$, it is not hard to see that $\PRk(\ModS(\vdash_{i})) \subseteq \Mod(\vdash_{i})$. In particular, this implies that $\SSS\PRk(\ModS(\vdash_{i})) \subseteq \Mod(\vdash_{i})$ and, therefore, that $\langle \A_{i}, F_{i} \rangle^{\ast} \in \RRR(\Mod(\vdash_{i}))$.

(ii)$\Rightarrow$(i): From the definition of $\bigotimes_{i \in I}{\vdash_{i}}$ it follows that $\bigotimes_{i \in I}\langle \A_{i}, F_{i} \rangle$ is a model of $\bigotimes_{i \in I}{\vdash_{i}}$. As submatrices of models are still models, this implies that $\langle \A, F \rangle \in \Mod(\bigotimes_{i \in I}{\vdash_{i}})$. Finally, the matrix $\langle \A, F \rangle$ is reduced by Corollary \ref{Cor:Taylor-Leibniz}.
\end{proof}

The following observation is well-known \cite[pag.\ 205]{AAL-AIT-f}.

\begin{Lemma}\label{Lem:no-theorems-trivial} Let $\vdash$ be a logic, and $\A$ an algebra.
\benroman
\item If $\langle \A, \emptyset \rangle \in \RRR(\Mod(\vdash))$, then $\A$ is the trivial algebra ${\bf 1}$.
\item A logic $\vdash$ has theorems if and only if $\langle \bf{1}, \emptyset \rangle \notin \ModS(\vdash)$ or, equivalently, if $\langle \bf{1}, \emptyset \rangle \notin \RRR(\Mod(\vdash))$.
\eroman
\end{Lemma}

As a consequence we obtain a transparent description of $\RRR(\Mod(\bigotimes_{i \in I}{\vdash_{i}}))$:

\begin{Proposition}
 Let $\{ \vdash_{i} \colon i \in I \}$ be a family of logics. The class $\RRR(\Mod(\bigotimes_{i \in I}{\vdash_{i}}))$ consists  of matrices satisfying condition (ii) of Lemma \ref{Lem:Taylor-reduced}, plus $\langle \boldsymbol{1}, \emptyset \rangle$ in case some $\vdash_{i}$ has no theorems.
\end{Proposition}

\begin{proof}
This is an easy consequence of Proposition \ref{Prop:Taylor-thms}, and of Lemmas \ref{Lem:Taylor-reduced} and \ref{Lem:no-theorems-trivial}
\end{proof}

Let $\{ \LL_{i} \colon i \in I \}$  be a family of languages and $\langle \A, F \rangle$ be a $\LL_{j}$-matrix for some $j \in I$. We denote by $\langle \A, F \rangle^{\flat}$ the $\bigotimes_{i \in I}\LL_{i}$-matrix $\bigotimes_{i \in I} \langle \A^{-}_{i}, F^{-}_{i}\rangle$, where
\[
\langle \A^{-}_{i}, F^{-}_{i}\rangle \coloneqq \left\{ \begin{array}{ll}
\langle \A, F \rangle & \text{if $i= j$}\\
  	\langle \boldsymbol{1}, \{ 1 \} \rangle& \text{otherwise.}\\
  \end{array} \right.
\]
Note that if $\langle \A, F \rangle$ is reduced, then  $\langle \A, F \rangle^{\flat}$ is reduced as well.

\begin{Lemma}\label{Lem:Taylor-Leibniz}
If $\{ \vdash_{i} \colon i \in I \}$ is a family of logics, 
\[
\RRR(\Mod(\bigotimes_{i \in I}{\vdash_{i}})) \subseteq \PSD(\bigotimes_{i \in I} \RRR(\Mod(\vdash_{i})))\subseteq \PSD\RRR(\Mod(\bigotimes_{i \in I}{\vdash_{i}})).
\]
\end{Lemma}

\begin{proof}
We detail only the proof of the first inclusion, since the proof of the second one exploits similar ideas. Consider a matrix $\langle \A, F \rangle \in \RRR(\Mod(\bigotimes_{i \in I}{\vdash_{i}}))$. First we consider the case where $F = \emptyset$. As the matrix $\langle \A, F \rangle$ is reduced, we know that $\A$ is trivial by Lemma \ref{Lem:no-theorems-trivial}(i). Now, the fact that $F$ is empty implies that $\bigotimes_{i \in I}{\vdash_{i}}$ has no theorems. From Proposition \ref{Prop:Taylor-thms} it follows that there is $j \in I$ such that $\vdash_{j}$ has no theorems. Therefore by Lemma \ref{Lem:no-theorems-trivial}(ii) the $\LL_{j}$-matrix $\langle \boldsymbol{1}, \emptyset \rangle$ belongs to $\RRR(\Mod( \vdash_{j}))$. As a consequence we obtain that
\[
\langle \A, F \rangle = \langle \boldsymbol{1}, \emptyset \rangle^{\flat} \in \bigotimes_{i \in I}(\RRR(\Mod(\vdash_{i})).
\] 

Then we consider the case where $F \ne \emptyset$. From Lemma \ref{Lem:Taylor-reduced} we know that $\langle \A, F \rangle \leq_{\textup{sd}} \bigotimes_{i \in I}\langle \A_{i}, F_{i} \rangle$ for some $\langle \A_{i}, F_{i}\rangle \in \RRR(\Mod(\vdash_{i}))$. Moreover, it is easy to see  that the map
\[
f \colon \prod_{i \in I} \langle \A_{i}, F_{i}\rangle^{\flat} \to \bigotimes_{i \in I}\langle \A_{i}, F_{i} \rangle
\]
defined by the rule
\[
f(\vec{a})(i) \coloneqq \vec{a}(i)(i)\text{, for every }i \in I,
\]
is an isomorphism. Together with the fact that $\langle \A, F \rangle \leq_{\textup{sd}} \bigotimes_{i \in I}\langle \A_{i}, F_{i} \rangle$, this implies that
\[
\langle \A, F \rangle \leq \prod_{i \in I} \langle \A_{i}, F_{i}\rangle^{\flat}
\]
is a subdirect product. Hence we conclude that $\langle \A, F \rangle \in\PSD(\bigotimes_{i \in I} \RRR(\Mod(\vdash_{i})))$.
\end{proof}


\begin{proof}[Proof of Proposition \ref{Prop:Taylor-Suszko}]
We begin by proving the first part. From Lemma \ref{Lem:PSD-Suszko} and \ref{Lem:Taylor-Leibniz} it follows that
\begin{align*}
\ModS(\bigotimes_{i \in I}{\vdash_{i}}) = \PSD \RRR(\Mod(\bigotimes_{i \in I}{\vdash_{i}})) \subseteq \PSD\PSD (\bigotimes_{i \in I} \RRR(\Mod(\vdash_{i}))) = \PSD (\bigotimes_{i \in I} \RRR(\Mod(\vdash_{i}))).
\end{align*}
Moreover, since $\RRR(\Mod(\vdash_{i})) \subseteq \ModS(\vdash_{i})$ for every $i \in I$, we have
\[
\PSD(\bigotimes_{i \in I} \RRR(\Mod(\vdash_{i})))  \subseteq  \PSD(\bigotimes_{i \in I} \ModS(\vdash_{i})).
\]

It only remains to prove that $\PSD(\bigotimes_{i \in I} \ModS(\vdash_{i})) \subseteq \ModS(\bigotimes_{i \in I}{\vdash_{i}})$. Since the class $\ModS(\bigotimes_{i \in I}{\vdash_{i}})$ is closed under subdirect products by Lemma \ref{Lem:PSD-Suszko}, it suffices to show that $\bigotimes_{i \in I} \ModS(\vdash_{i}) \subseteq \ModS(\bigotimes_{i \in I}{\vdash_{i}})$. To this end, consider a matrix $\langle \A_{i}, F_{i} \rangle  \in \ModS(\vdash_{i})$ for each $i \in I$. By Lemma \ref{Lem:PSD-Suszko},  for every $i \in I$ there is a family $\{ \langle \A_{i}^{j}, F_{i}^{j}\rangle \colon j \in J_{i} \} \subseteq \RRR(\Mod(\vdash_{i}))$ such that $\langle \A_{i}, F_{i} \rangle \leq \prod_{j \in J_{i}}\langle \A_{i}^{j}, F_{i}^{j}\rangle$ is a subdirect product. We can assume without loss of generality that $J_{i} = J_{j}$ for every $i, j \in I$ (for instance, by adding trivial matrices to the factors of products when necessary). Accordingly, we drop the index $i$ in each $J_{i}$, and write simply $J$. Under this convention, it is easy to see that
\[
\bigotimes_{i \in I}\langle \A_{i}, F_{i} \rangle \leq \prod_{j \in J} (\bigotimes_{i \in I} \langle \A_{i}^{j}, F_{i}^{j}\rangle)
\]
is a subdirect product. Together with Lemmas \ref{Lem:Taylor-Leibniz} and \ref{Lem:PSD-Suszko}, this yields
\[
\bigotimes_{i \in I}\langle \A_{i}, F_{i} \rangle \in \PSD(\bigotimes_{i \in I}\RRR(\Mod(\vdash_{i}))) \subseteq \PSD\RRR(\Mod(\bigotimes_{i \in I}{\vdash_{i}})) = \ModS(\bigotimes_{i \in I}{\vdash_{i}}).
\]
Hence we conclude that $\bigotimes_{i \in I} \ModS(\vdash_{i}) \subseteq \ModS(\bigotimes_{i \in I}{\vdash_{i}})$, as desired.

To prove the second part,  we rely on the first one. Consider a matrix $\langle \A, F \rangle \in \ModS(\bigotimes_{i \in I}{\vdash_{i}}) = \PSD(\bigotimes_{i \in I} \ModS(\vdash_{i}))$. We can assume without loss of generality that
\[
\langle \A, F \rangle \subseteq \prod_{j \in J} \bigotimes_{i \in I} \langle \B_{i}^{j}, G_{i}^{j}\rangle
\]
is a subdirect product for some families $\{\langle \B_{i}^{j}, G_{i}^{j}\rangle \colon j \in J\} \subseteq \ModS(\vdash_{i})$, one for each $i \in I$. Then for every $i \in I$, let
\begin{align*}
f \colon \prod_{j \in J} (\bigotimes_{i \in I} \langle \B_{i}^{j}, G_{i}^{j}\rangle) &\to  \bigotimes_{i \in I} (\prod_{j \in J}\langle \B_{i}^{j}, G_{i}^{j}\rangle)\\
\pi_i \colon \prod_{i \in I}\prod_{j \in J}B_{i}^{j} &\to \prod_{j \in J}B_{i}^{j}
\end{align*}
be, respectively, the isomorphism defined in the proof of Lemma \ref{Lem:Taylor-reduced}, and the natural projection on the $i$-th component. Bearing this in mind, for every $i \in I$ let $\langle \C_{i}, H_{i}\rangle$ be the matrix where $\C_i$ is the subalgebra $\pi_i[f[\A]] \subseteq \prod_{j \in J}\B_{i}^{j}$ and $H_i = \pi_i[f[F]]$. 

The restriction $f{\upharpoonright}_A \colon \langle \A, F \rangle \to \bigotimes_{i \in I} \langle \C_{i}, H_{i}\rangle$ is a well-defined matrix embedding such that $\pi_i[f{\upharpoonright}_A [A]] = C_i$ for every $i \in I$. Hence, we conclude that
\[
\langle \A, F\rangle \leq_{\textup{sd}} \bigotimes_{i \in I} \langle \C_{i}, H_{i}\rangle. 
\]
Now, it is not hard to see that $\langle \C_{i}, H_{i}\rangle$ is a subdirect product of $\prod_{j \in J}\langle \B_i^{j}, G_{i}^{j}\rangle$, for every $i \in I$. Since each $\ModS(\vdash_{i})$ is closed under subdirect products, this implies that $\langle \C_{i}, H_{i}\rangle \in \ModS(\vdash_i)$ for every $i \in I$. This proves the inclusion from left to right.

To prove the other inclusion, let $\langle \A, F\rangle \leq_{\textup{sd}}  \bigotimes_{i \in I}\langle \A_{i}, F_{i}\rangle$ where $\langle \A_{i}, F_{i}\rangle\in  \ModS(\vdash_i)$ for every $i \in I$. Then, as in the proof of Lemma \ref{Lem:Taylor-Leibniz}, we have that
$\langle \A, F \rangle$ is a subdirect product of  $\prod_{i \in I} \langle \A_{i}, F_{i}\rangle^{\flat}$. From the fact that $\langle \A_{i}, F_{i}\rangle^{\flat} \in \bigotimes_{i \in I} \ModS(\vdash_{i})$ for every $i \in I$, we obtain that $\langle \A, F \rangle \in \PSD(\bigotimes_{i \in I} \ModS(\vdash_{i}))$.
\end{proof}

The characterization of $\ModS(\bigotimes_{i \in I}{\vdash_{i}})$ given in Proposition \ref{Prop:Taylor-Suszko} has a particularly appealing simplification in the case where the index set $I$ is finite.

\begin{Corollary}\label{Cor:FiniteCase}
If $\vdash$ and $\vdash'$ are logics, then
\[
\ModS({\vdash} \bigotimes {\vdash'}) = \ModS(\vdash) \bigotimes \ModS(\vdash').
\]
\end{Corollary}

\begin{proof}
As shown essentially in \cite[Lem.\ 1.9 and 1.10]{Tay73}, if $\class{K}_{1}$ and $\class{K}_{2}$ are classes of matrices (resp.\ algebras) closed under subdirect products, then so is $\class{K}_{1} \bigotimes \class{K}_{2}$. Together with Lemma \ref{Lem:PSD-Suszko}, this implies that the class $\ModS(\vdash) \bigotimes \ModS(\vdash')$ is closed under subdirect products. By Proposition \ref{Prop:Taylor-Suszko} we conclude that $\ModS({\vdash} \bigotimes {\vdash'}) = \ModS(\vdash) \bigotimes \ModS(\vdash')$.
\end{proof}

\section{Finite suprema need not exist}

It is well known that if $\mathbb{A}$ is a poset whose universe is a set and in which infima of sets exist, then $\mathbb{A}$ is a complete lattice. Unfortunately, the proof of this fact   relies on the assumption that the universe of $\mathbb{A}$ is a set and, therefore, cannot be applied to the poset $\class{Log}$ (which is known to have infima of sets by Theorem \ref{Thm:Infima}). The situation is entirely different for $\class{Log}$. This  section is devoted to prove the following:

\begin{Theorem}\label{Thm:NotJoins}
Finite suprema need not exist in $\class{Log}$.
\end{Theorem}

The proof of this theorem builds on a counterexample.\ Let $\A = \langle A; \lor, \bold{a}, \bold{b}, \bold{0}\rangle$ be the join-semilattice, expanded with constants\footnote{ We use in this section ``constant(s)'' as an abbreviation for ``constant unary operation(s)''.}, depicted below:
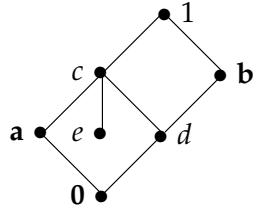
\[
\xymatrix@R=23pt @C=23pt @!0{
&& *-{ \  \ \ \bullet  \ 1 } \ar@{-}[dr]\ar@{-}[dl]&\\
& *-{c \ \bullet  \ \ \ }\ar@{-}[dr]\ar@{-}[dl] && *-{ \?  \ \ \bullet  \ \bold{b}}\ar@{-}[dl]\\
*-{\bold{a} \ \bullet  \ \ \ } \ar@{-}[dr]&*-{e \ \bullet  \ \ \ }\ar@{-}[u]& *-{\?  \ \ \bullet  \ d } \ar@{-}[dl]&\\
&*-{ \bold{0} \ \bullet  \ \ \ }&&
}
\]
Then let $\vdash_{\lor}$ be the logic in countably many variables induced by the set of matrices
\[
\{ \langle \A, \{ 1 \} \rangle, \langle \A, \{ 1, c \} \rangle \}.
\]

\begin{Fact}\label{Fact:TwoMatrices}
We have that $\langle \A, \{ 1 \}\rangle \in \ModS(\vdash_{\lor})$.
\end{Fact}

\begin{proof}
It is clear that $\langle \A, \{ 1 \}\rangle$ is a model of $\vdash_{\lor}$. Hence it will be enough to prove that $\tarski_{\vdash_{\lor}}^{\A}\{ 1 \}$ is the identity relation on $A$. From the definition of $\vdash_{\lor}$ it follows that $\{ c, 1\}$ is a deductive filter of $\vdash_{\lor}$ on $\A$. Now, an easy computation shows that:
\benroman
\item The blocks of $\leibniz^{\A}\{ 1 \}$ are $\{ a, e, c \}, \{ 0, d\}, \{ b \}, \{ 1\}$.
\item The blocks of $\leibniz^{\A}\{ c, 1 \}$ are $\{ 0\}, \{ a\}, \{e \}, \{ b, d\}, \{ c, 1\}$.
\eroman
Together with the fact that
\[
\tarski_{\vdash_{\lor}}^{\A}\{ 1 \} \subseteq \leibniz^{\A}\{ 1 \} \cap \leibniz^{\A}\{ c, 1 \},
\]
this implies that $\tarski_{\vdash_{\lor}}^{\A}\{ 1 \}$ is the identity relation on $A$.
\end{proof}

\begin{Fact}\label{Fact:four}
The algebraic reducts  of the matrices in $\ModS(\vdash_{\lor})$ are either trivial or have at least four elements.
\end{Fact}

\begin{proof}
In this proof we assume that semilattices are equipped with the join-order. Consider a matrix $\langle \B, F\rangle \in \ModS(\vdash_{\lor})$ such that $\B$ is non-trivial. By Corollary \ref{Cor:Alg-general-trick} we know that $\B$ is a semilattice with constants $\bold{a}, \bold{b}, \bold{0}$ such that
\begin{equation}\label{Eq:Below}
\bold{0} \leq \bold{a} \text{ and }\bold{0} \leq  \bold{b}.
\end{equation}
Now, since $\B$ is non-trivial, we know that $F \ne B$. Together with the fact that
\[
\bold{a} \vdash_{\lor} x \quad \bold{b} \vdash_{\lor} x \quad \bold{0} \vdash_{\lor} x
\]
this implies that $\bold{a}, \bold{b}, \bold{0} \notin F$. Observe that $\bold{a} \lor \bold{b} \in F$, since $\emptyset \vdash_{\lor}\bold{a} \lor \bold{b}$. Hence, to conclude that $B$ has at least four elements, it will be enough to check that $\bold{a}, \bold{b}, \bold{0}$ are different one from the other. From the fact that $\bold{a}, \bold{b} \notin F$ and $\bold{a} \lor \bold{b} \in F$, it follows that $\bold{a}$ and $\bold{b}$ are incomparable in the order of $\B$. Together with (\ref{Eq:Below}), this implies that $\bold{0}$ is different from $\bold{a}$ and $\bold{b}$.
\end{proof}

We say that a \textit{negation algebra} is an algebra $\B = \langle B; \lnot\rangle$ where $\lnot$ is a unary operation with at most one fix point, and such that $\lnot \lnot a = a $ for all $a \in B$. We denote by $\class{NA}$ the class of negation algebras, and by $\vdash_{\lnot}$ be the negation fragment of classical propositional logic (formulated in countably many variables). The relation between $\vdash_{\lnot}$ and $\class{NA}$ is captured by the following result:

\begin{Fact}\label{Fact:ModNeg}
$\ModS(\vdash_{\lnot})$ is the class of matrices $\langle \B, F\rangle$ such that either $\B$ is trivial or $\B$ is a negation algebra and in this case  either $F= \emptyset$ or $F= \{ a \}$ for some $a \in B$ that is not a fixed point of $\lnot$.
\end{Fact}

\begin{proof}
The interested reader may consult the Appendix for the details.
\end{proof}


Now, given a cardinal $\kappa > 0$ and $\alpha < \kappa$, we let $\A_{\alpha, \kappa}$ be the expansion of $\A$ with a  constant   for every element of $A$, a unary operation $\lnot$ defined as
\[
\lnot 0 = d  \quad \lnot d = 0 \quad \lnot a = c  \quad \lnot c = a \quad \lnot 1 = b  \quad \lnot b = 1 \quad \lnot e = e,
\]
and with a set of binary operations $\{ \multimap_{\beta} \colon \beta < \kappa \}$ defined for every $\beta < \kappa$ and $p, q \in A$ as follows:
\begin{displaymath}
p \multimap_{\beta}q\coloneqq \left\{\begin{array}{@{\,}ll}
1 & \text{if $p=q$ or $\beta \ne \alpha$}\\
0 & \text{if $p \ne q$ and $\beta = \alpha$.}\\
\end{array} \right.
\end{displaymath}
Then let $\vdash_{\kappa}$ be the logic formulated in countably many variables induced by the class of matrices $\{ \langle \A_{\alpha, \kappa}, \{ 1 \} \rangle \colon \alpha < \kappa \}$.

\begin{Fact}\label{Fact:Algebraizable}
For every $\kappa > 0$, the logic $\vdash_{\kappa}$ is equivalential.
\end{Fact}
\begin{proof}
Consider the set
\[
\Delta(x, y) \coloneqq \{ x \multimap_{\alpha} y \colon \alpha < \kappa \}.
\]
It is easy to see that $\Delta$ witnesses the validity of the rules in Theorem \ref{Thm:congruence-formulas}. Hence we conclude that $\vdash_{\kappa}$ is equivalential.
\end{proof}

%

\begin{Fact}\label{Fact2} For every $\kappa > 0$, $\llbracket \vdash_{\kappa}\rrbracket$ is an upper bound of $\llbracket\vdash_{\lor}\rrbracket$ and $\llbracket\vdash_{\lnot}\rrbracket$ in $\class{Log}$.
\end{Fact}

\begin{proof}
Consider the class of matrices $\class{K} \coloneqq \{ \langle \A_{\alpha, \kappa}, \{ 1 \} \rangle \colon \alpha < \kappa \}$. It is clear that $\vdash_{\kappa}$ is the logic induced by $\class{K}$ and that $\SSS(\class{K}) = \class{K}$. Then let $\btau$ be the identity translation of $\LL_{\vdash_{\lor}}$ into $\LL_{\vdash_{\kappa}}$. By Fact \ref{Fact:TwoMatrices} we have that
\[
\langle \B^{\btau}, F \rangle = \langle \A, \{ 1 \} \rangle \in \ModS(\vdash_{\lor})
\]
for every $\langle \B, F \rangle \in \class{K}$.  Together with Fact \ref{Fact:Algebraizable} and Proposition \ref{Prop:alg-translations}, this implies that $\btau$ is an interpretation of $\vdash_{\lor}$ into $\vdash_{\kappa}$.

A similar argument (requiring Fact \ref{Fact:ModNeg}) shows that $\vdash_{\lor}$ is also interpretable in $\vdash_{\kappa}$.
\end{proof}

Suppose, with a view to contradiction, that there exists the supremum of $\llbracket \vdash_{\lor} \rrbracket$ and $\llbracket \vdash_{\lnot} \rrbracket$ in $\class{Log}$, i.e. that there exists a logic $\vdash$ such that 
\begin{equation}\label{Eq:Counterexample}
\llbracket \vdash \rrbracket = \llbracket \vdash_{\lor} \rrbracket \lor \llbracket \vdash_{\lnot} \rrbracket.
\end{equation}
From now on, our aim is to obtain a contradiction.

\begin{Fact}\label{Fact:Kappa}
For every $\kappa>0$, we have $\llbracket \vdash \rrbracket \leq \llbracket \vdash_{\kappa}\rrbracket$.
\end{Fact}

\begin{proof}
This is a direct consequence of Fact \ref{Fact2}.
\end{proof}

Now, since $\vdash$ is a logic, its language is a \textit{set}, say of cardinality $\kappa$. We can assume without loss of generality that $\kappa$ is infinite (if it is not, then we can add to it infinitely many unary operations whose interpretation in $\ModS(\vdash)$ would be the identity map).\ By Fact \ref{Fact:Kappa} there is an interpretation $\btau$ of $\vdash$ into $\vdash_{\kappa^{+}}$.

\begin{Fact}\label{Eq:Multimap}
There is $\alpha < \kappa^{+}$ such that the symbol $\multimap_{\alpha}$ does not appear in the terms $\{ \btau(\varphi) \colon \varphi \in \LL_{\vdash} \}$.
\end{Fact}

\begin{proof}
Straightforward.
\end{proof}

From now on we will work with the special $\alpha < \kappa^{+}$ provided by Fact \ref{Eq:Multimap}. 
Let $\hat{\A}_{\alpha, \kappa^{+}}$ be the $\{ \multimap_{\alpha} \}$-free reduct of $\A_{\alpha, \kappa^{+}}$.  

\begin{Fact}\label{Fact:term-equivalence}
The algebra $\hat{\A}_{\alpha, \kappa^{+}}$ is term-equivalent to
\[
\langle A; \lor^{\A_{\alpha, \kappa^{+}}}, \lnot^{\A_{\alpha, \kappa^{+}}}, \bold{a}, \bold{b}, \bold{0}, e\rangle.
\]
\end{Fact}

\begin{proof}
Using negation, it is easy to see that all constants from  $\hat{\A}_{\alpha, \kappa^{+}}$ are definable in the displayed algebra. Moreover, if $\beta \ne \alpha$, then $\multimap_{\beta}^{\A_{\alpha, \kappa^{+}}}$ is a  constant map. This shows that all term-functions of $\hat{\A}_{\alpha, \kappa^{+}}$ are also term-functions of the algebra in the display. The converse is obvious.
\end{proof}

In what follows we will work under the identification of $\hat{\A}_{\alpha, \kappa^{+}}$ with the algebra displayed in Fact \ref{Fact:term-equivalence}. 

\begin{Fact}\label{Fact:no-constants}
If $\gamma(x)$ is a formula of $\hat{\A}_{\alpha, \kappa^{+}}$ such that $\langle A; \gamma^{\hat{\A}_{\alpha, \kappa^{+}}} \rangle$ is a negation algebra, then $\gamma$ can be obtained as a composition of $\lor$ and  $\lnot$.
\end{Fact}

\begin{proof}
Assume that  $\langle A; \gamma^{\hat{\A}_{\alpha, \kappa^{+}}} \rangle$ is a negation algebra.  Suppose, with a view to contradiction, that either $\bold{0}$ or $\bold{a}$ or $\bold{b}$ occur in $\gamma$. It is not hard to see that this implies that $e \notin \gamma^{\hat{\A}_{\alpha, \kappa^{+}}}[A]$. However, since $\langle A; \gamma^{\hat{\A}_{\alpha, \kappa^{+}}}\rangle$ is a negation algebra, we know that
\[
e = \gamma^{\hat{\A}_{\alpha, \kappa^{+}}} \gamma^{\hat{\A}_{\alpha, \kappa^{+}}}(e) \in \gamma^{\hat{\A}_{\alpha, \kappa^{+}}}[A],
\]
which is false. Hence we conclude that $\bold{0}$, $\bold{a}$, and $\bold{b}$ do not occur in $\gamma$. 

It only remains to prove that $e$ does not occur in $\gamma$. Suppose the contrary, with a view to contradiction. An easy induction on the construction  of formulas shows that if $\varphi(x)$ is a formula of $\hat{\A}_{\alpha, \kappa^{+}}$ in which $\bold{0}$, $\bold{a}$, and $\bold{b}$ do not occur and in which $e$ occurs, then $\varphi^{\hat{\A}_{\alpha, \kappa^{+}}}(0) \in \{ e, \bold{a}, c \}$. As $\bold{0}$, $\bold{a}$, and $\bold{b}$ do not occur in the composition $\gamma(\gamma(x))$, this means that $\gamma^{\hat{\A}_{\alpha, \kappa^{+}}}\gamma^{\hat{\A}_{\alpha, \kappa^{+}}}(0) \ne 0$. But this contradicts the fact that $\langle A; \gamma^{\hat{\A}_{\alpha, \kappa^{+}}}\rangle$ is a negation algebra, as desired.
\end{proof}

\begin{Fact}\label{Fact:LeibnizCongruence}
The blocks of $\leibniz^{\A_{\alpha, \kappa^{+}}^{\btau}}\{ 1 \}$ are $\{ 0, d \}, \{a, c, e \}, \{ b \}, \{ 1 \}$.
\end{Fact}

\begin{proof}
By Fact \ref{Eq:Multimap} we know that the term-functions of $\A^{\btau}_{\alpha, \kappa^{+}}$ are also term-functions of $\hat{\A}_{\alpha, \kappa^{+}}$. In particular, this means that $\Con \hat{\A}_{\alpha, \kappa^{+}} \subseteq \Con \A_{\alpha, \kappa^{+}}^{\btau}$. 

Consider the equivalence relation $\theta$ on $A$ determined by the partition in the statement. Using for instance Fact \ref{Fact:term-equivalence}, it is easy to see that $\theta$ is a congruence of $\hat{\A}_{\alpha, \kappa^{+}}$. Then $\theta$ is also a congruence of $\A^{\btau}_{\alpha, \kappa^{+}}$. As $\theta$ is compatible with $\{ 1 \}$, this implies that $\theta \subseteq \leibniz^{\A_{\alpha, \kappa^{+}}^{\btau}}\{ 1 \}$.

As a consequence, we obtain that $A / \leibniz^{\A_{\alpha, \kappa^{+}}^{\btau}}\{ 1 \}$ is a set of at most four elements. Moreover, since $\leibniz^{\A_{\alpha, \kappa^{+}}^{\btau}}\{ 1 \}$ is compatible with $\{ 1\}$, we know that $\langle 0, 1 \rangle \notin \leibniz^{\A_{\alpha, \kappa^{+}}^{\btau}}\{ 1 \}$. Therefore we have
\begin{equation}\label{Eq:cardinality}
2 \leq \vert A / \leibniz^{\A_{\alpha, \kappa^{+}}^{\btau}}\{ 1 \} \vert \leq 4.
\end{equation}
Now, it is easy to see that $\langle \A_{\alpha, \kappa^{+}}, \{ 1 \} \rangle \in \RRR(\ModS(\vdash_{\kappa^{+}}))$. Since $\btau$ is an interpretation of $\vdash$ into $\vdash_{\kappa^{+}}$, this implies that $\langle \A_{\alpha, \kappa^{+}}^{\btau}, \{ 1 \} \rangle \in \Mod(\vdash)$ and, therefore, that $\langle \A_{\alpha, \kappa^{+}}^{\btau}, \{ 1 \} \rangle^{\ast} \in \ModS(\vdash)$. Together with Fact \ref{Fact:four} and ${\vdash_{\lor}} \leq {\vdash}$, this implies that either the matrix $\langle \A_{\alpha, \kappa^{+}}^{\btau}, \{ 1 \} \rangle^{\ast}$ is trivial, or the congruence $\leibniz^{\A_{\alpha, \kappa^{+}}^{\btau}}\{ 1 \}$ has at least four blocks. By (\ref{Eq:cardinality}) we conclude that $\leibniz^{\A_{\alpha, \kappa^{+}}^{\btau}}\{ 1 \}$ has exactly four blocks. Together with the fact that $\theta \subseteq \leibniz^{\A_{\alpha, \kappa^{+}}^{\btau}}\{ 1 \}$, this implies that $\theta = \leibniz^{\A_{\alpha, \kappa^{+}}^{\btau}}\{ 1 \}$.
\end{proof}

We are now ready to produce the desired contradiction. To this end, recall that there is an interpretation $\brho$ of $\vdash_{\lnot}$ into $\vdash$. Since $\langle \A_{\alpha, \kappa^{+}}, \{ 1 \} \rangle \in \ModS(\vdash_{\kappa^{+}})$, we can apply Fact \ref{Fact:ModNeg} obtaining that $\langle A; \btau\brho(\lnot)^{\A_{\alpha, \kappa^{+}}} \rangle$ is a negation algebra. By Fact \ref{Eq:Multimap} we know that the function $\btau\brho(\lnot)^{\A_{\alpha, \kappa^{+}}}\colon A \to A$ is a term-function of $\hat{\A}_{\alpha, \kappa^{+}}$. Hence we can apply Fact \ref{Fact:no-constants} obtaining that $\btau\brho(\lnot)^{\A_{\alpha, \kappa^{+}}}$ can be produced as a composition of the functions 
\[
\lnot^{\A_{\alpha, \kappa^{+}}} \colon A \to A \text{ and }\lor^{\A_{\alpha, \kappa^{+}}} \colon A \times A \to A.
\] 
This yields that
\begin{equation}\label{Eq:contradiction}
\btau\brho(\lnot)^{\A_{\alpha, \kappa^{+}}}(0) \in \{ 0, d\} \text{ and }\btau\brho(\lnot)^{\A_{\alpha, \kappa^{+}}}(e)  = e.
\end{equation}

From the fact that $\langle \A_{\alpha, \kappa^{+}}, \{ 1 \} \rangle \in \ModS(\vdash_{\kappa^{+}})$ it follows that $\langle \A_{\alpha, \kappa^{+}}^{\btau}, \{ 1 \}\rangle \in \Mod(\vdash)$. In particular, this implies that 
\[
\langle \A_{\alpha, \kappa^{+}}^{\btau}/ \theta, \{ 1 \} / \theta \rangle \in \ModS(\vdash),
\]
where $\theta \coloneqq \leibniz^{\A_{\alpha, \kappa^{+}}^{\btau}}\{ 1 \}$. Together with Fact \ref{Fact:ModNeg}, this yields that $\langle A/ \theta; \btau\brho(\lnot)^{\A_{\alpha, \kappa^{+}}/ \theta}\rangle$ is a negation algebra. However, by Fact \ref{Fact:LeibnizCongruence} and (\ref{Eq:contradiction}) this negation algebra has two distinct fixed points for negation (namely $0/\theta$ and $e/ \theta$), which is impossible. Hence we reached a contradiction, establishing Theorem \ref{Thm:NotJoins}.


\section{The lattice of equivalential logics}

Even if  suprema need not exist in $\class{Log}$ there is an important subsemilattice of $\class{Log}$ where suprema exist, i.e.\ the lattice of equivalential logics. 

\begin{Proposition}\label{Prop:Equiv-interpretation} \
\benroman
\item Let $\vdash$ and $\vdash'$ be logics. If $\vdash$ is equivalential and ${\vdash} \leq {\vdash'}$, then $\vdash'$ is also equivalential.
\item If $\{ \vdash_{i} \colon i \in I \}$ is a family of equivalential logics, then $\bigotimes_{i \in I} {\vdash_{i}}$ is equivalential.
\eroman
\end{Proposition}

\begin{proof}
(i): Let $\Delta(x, y)$ be the set of formulas witnessing the fact that $\vdash$ is equivalential, as in Theorem \ref{Thm:congruence-formulas}. Moreover, let $\btau$ be an interpretation of $\vdash$ into $\vdash'$. We consider the set $\Sigma(x, y) \coloneqq \btau[\Delta]$ of formulas of $\LL_{\vdash'}$. In order to establish that $\vdash'$ is equivalential, it will be enough to show that $\Sigma$ and $\vdash'$ satisfy the conditions in Theorem \ref{Thm:congruence-formulas}.

From Proposition \ref{Prop:interpretation} it follows that $\emptyset \vdash' \Sigma(x, x)$ and $x, \Sigma(x, y) \vdash' y$. It only remains to prove that for every $n$-ary connective $\ast$ of $\vdash'$, 
\begin{equation}\label{Eq:Equiv-syntactic-proof}
\bigcup_{1 \leq i \leq n}\Sigma(x_{i}, y_{i}) \vdash' \Sigma(\ast(x_{1}, \dots, x_{n}), \ast(y_{1}, \dots, y_{n})).
\end{equation}
To this end, consider an $n$-ary connective $\ast$ of $\vdash'$, a matrix $\langle \A, F \rangle \in \ModS(\vdash')$, and tuples $\vec{a}, \vec{c} \in A^{n}$ such that
\[
\bigcup_{1 \leq i \leq n}\Sigma^{\A}(a_{i}, c_{i}) \subseteq F.
\]
Since $\Sigma = \btau[\Delta]$, we have 
\[
\bigcup_{1 \leq i \leq n}\Delta^{\A^{\btau}}(a_{i}, c_{i}) \subseteq F.
\]
As $\Delta$ is a set of congruence formulas for $\vdash$, and $\langle \A^{\btau}, F \rangle \in \ModS(\vdash) = \RRR(\Mod(\vdash))$, the above display implies that $\vec{a} = \vec{c}$. As a consequence, we obtain that $\ast^{\A}(\vec{a}) = \ast^{\A}(\vec{c}\?)$. Since $\emptyset \vdash' \Sigma(x, x)$ and $\langle \A, F \rangle \in \Mod(\vdash')$, this yields 
\[
\Sigma^{\A}(\ast(\vec{a}), \ast(\vec{c}\?)) = \Sigma^{\A}(\ast(\vec{a}), \ast(\vec{a})) \subseteq F.
\]
Hence we conclude that (\ref{Eq:Equiv-syntactic-proof}) holds.

(ii): Given $i \in I$, let $\Delta_{i}(x, y)$ be a set of congruence formulas for $\vdash_{i}$. Observe that the Cartesian product $\prod_{i \in I} \Delta_{i}$ can be viewed as a set $\Delta(x, y)$ of formulas of $\bigotimes_{i \in I}{\vdash_{i}}$. Since the various $\Delta_{i}$ satisfy the rules in Theorem \ref{Thm:congruence-formulas}, and $\bigotimes_{i \in I}{\vdash_{i}}$ is the logic induced by $\bigotimes_{i \in I} \ModS(\vdash_{i})$, it is easy to see that the set $\Delta$ satisfies the rules in Theorem \ref{Thm:congruence-formulas} as well. As a consequence, we conclude that $\bigotimes_{i \in I} {\vdash_{i}}$ is an equivalential logic.
\end{proof}

The above result motivates the following definition:

\begin{law}
Let $\class{Equiv}$ be the subposet of $\class{Log}$ that contains the classes $\llbracket \vdash \rrbracket$ such that $\vdash$ is an equivalential logic.
\end{law}

From Proposition \ref{Prop:Equiv-interpretation} it follows that $\class{Equiv}$ is a \textit{set-complete filter} of $\class{Log}$, i.e.\ an upset that is closed under infima of sets. Moreover, we shall prove that in $\class{Equiv}$ suprema of sets exist.

\begin{law}
Given a family $\{\LL_i \colon i \in I \}$ of languages, we let $\bigoplus_{i \in I}{\LL_i}$ be the language consisting of the disjoint union of the various $\LL_{i}$. Moreover, given a family $\{ \vdash_{i} \colon i \in I \}$ of equivalential logics, we let $\bigoplus_{i \in I}{\vdash_{i}}$ be the logic in the language $\bigoplus_{i \in I}{\LL_i}$ formulated in $\Sigma_{i \in I} \kappa_{\vdash_{i}}$ variables and induced by the following class of $\bigoplus_{i \in I}{\LL_{\vdash_{i}}}$-matrices:
\begin{equation}\label{Eq:Supremum}
  \{ \langle \A, F\rangle \colon  \text{the $\LL_{\vdash_{i}}$-reduct of $\langle \A, F\rangle$ belongs to $\RRR(\Mod(\vdash_{i}))$ for all $i \in I$}	\}.
\end{equation}
\end{law}
We will show that $\llbracket \bigoplus_{i \in I}{\vdash_{i}}\rrbracket$ is the supremum of $\{ \llbracket \vdash_{i}\rrbracket \colon i \in I \}$ both in $\class{Log}$ and $\class{Equiv}$.

%
%

\begin{Lemma}\label{Lem:Equiv-tricks}
Let $\{ \vdash_{i} \colon i \in I \}$ be a family of equivalential logics.
\benroman 
\item If $\Delta$ is a set of congruence formulas for $\vdash_{i}$, then so it is for $\bigoplus_{i \in I}{\vdash_{i}}$.
\item The logic $\bigoplus_{i \in I}{\vdash_{i}}$ is equivalential.
\item $\ModS(\bigoplus_{i \in I}{\vdash_{i}})$ is the class of matrices in (\ref{Eq:Supremum}).
\item $\llbracket \bigoplus_{i \in I}{\vdash_{i}} \rrbracket$ is the supremum of $\{ \llbracket \vdash_{i}\rrbracket \colon i \in I\}$ both in $\class{Equiv}$ and in  $\class{Log}$.
\eroman
\end{Lemma}

\begin{proof}
(i): Observe that the $\LL_{\vdash_{i}}$-reducts of the matrices in (\ref{Eq:Supremum}) are reduced. Together with the fact that $\Delta$ is a set of congruence formulas for $\vdash_{i}$, this easily implies that $\Delta$ satisfies the conditions of Theorem \ref{Thm:congruence-formulas} for $\bigoplus_{i \in I}{\vdash_{i}}$. As a consequence, we conclude that $\Delta$ is a set of congruence formulas for $\bigoplus_{i \in I}{\vdash_{i}}$.

(ii): Immediate from (i). (iii): Let $\class{M}$ be the class of matrices in (\ref{Eq:Supremum}). It is easy to see that the matrices in $\class{M}$ are reduced and, therefore, that $\class{M} \subseteq \ModS(\bigoplus_{i \in I}{\vdash_{i}})$. To prove the other inclusion, consider $\langle \A, F \rangle \in \ModS(\bigoplus_{i \in I}{\vdash_{i}})$. As $\bigoplus_{i \in I}{\vdash_{i}}$ is equivalential by (ii), we can apply Theorem \ref{Thm:congruence-formulas} obtaining that $\langle \A, F \rangle \in \RRR(\Mod(\bigoplus_{i \in I}{\vdash_{i}}))$. It will be enough to show that (for every $i \in I$) the $\LL_{\vdash_{i}}$-reduct $\langle \A^{-}, F \rangle$ of $\langle \A, F \rangle$ is a reduced model of $\vdash_{i}$. The fact that $\langle \A^{-}, F\rangle$ is a model of $\vdash_{i}$ is clear. To prove that it is reduced, let $\Delta$ be a set of congruence formulas of $\vdash_{i}$. By (i) we know that $\Delta$ is also a set of congruence formulas for $\bigoplus_{i \in I}{\vdash_{i}}$. Together with the fact that $\langle \A, F \rangle$ is a reduced model of $\bigoplus_{i \in I}{\vdash_{i}}$, this implies that for every $a, b \in A$,
\[
a = b \Longleftrightarrow \Delta^{\A}(a, b) \subseteq F \Longleftrightarrow \Delta^{\A^{-}}(a, b) \subseteq F.
\]
Since $\langle \A^{-}, F \rangle$ is a model of $\vdash_{i}$, this implies that the matrix $\langle \A^{-}, F \rangle$ is reduced.

(iv): By (i) we know that $\llbracket \bigoplus_{i \in I}{\vdash_{i}} \rrbracket$ belongs to $\class{Equiv}$. Hence it will be enough to show that it is the supremum of $\{ \llbracket \vdash_{i}\rrbracket \colon i \in I\}$ in $\class{Log}$. Recall from Theorem \ref{Thm:congruence-formulas} that $\ModS(\vdash_{i}) = \RRR(\Mod(\vdash_{i}))$ for all $i \in I$. Together with (iii), this implies that ${\vdash_{j}} \leq {\bigoplus_{i \in I}{\vdash_{i}}}$ for all $j \in I$.

Then consider a logic $\vdash$ such that ${\vdash_{i}} \leq {\vdash}$ for every $i \in I$. Then for every $i \in I$, there is an interpretation $\btau_{i}$ of $\vdash_{i}$ into $\vdash$. Observe that all these $\btau_{i}$ can be joined together into a translation $\btau$ of $\bigoplus_{i \in I}{\LL_{i}}$ into $\LL_{\vdash}$. We will show that $\btau$ is also an interpretation of $\bigoplus_{i \in I}{\vdash_{i}}$ into $\vdash$. To this end, consider a matrix $\langle \A, F\rangle \in \ModS(\vdash)$. We know that $\langle \A^{\btau_{i}}, F\rangle \in \ModS(\vdash_{i}) = \RRR(\Mod(\vdash_{i}))$ for every $i \in I$. This implies that the matrix $\langle \A^{\btau}, F\rangle$ belongs to the class in (\ref{Eq:Supremum}). By (iii) we conclude that $\langle \A, F\rangle \in \ModS(\bigoplus_{i \in I}{\vdash_{i}})$.
\end{proof}

As a consequence, we obtain the following:

\begin{Theorem}\label{Thm:ProtoalgebraicLattice}
$\class{Equiv}$ is a set-complete lattice, i.e.\ infima and suprema of subsets of $\class{Equiv}$ exist. Moreover, these infima and suprema coincide with those of $\class{Log}$.
\end{Theorem}

\begin{proof}
From Proposition \ref{Prop:Equiv-interpretation}(ii) and Lemma \ref{Lem:Equiv-tricks}(iv).
\end{proof}

\begin{problem}
Do suprema of protoalgebraic logics \cite{AAL-AIT-f} exist as well?
\end{problem}

An adaptation of an argument given in \cite[pag. 34]{GaTa84} shows that the lattice $\class{Equiv}$ is not modular. However, to our knowledge, the following problem remains open:

\begin{problem}
Do $\class{Equiv}$ and $\class{Var}$ satisfy any non-trivial lattice equation?
\end{problem}

\section{The top and the bottom}

In this section we will describe the top and the bottom parts of $\class{Log}$. To this end, recall that a logic $\vdash$ is \textit{inconsistent} if $\Gamma \vdash \varphi$ for every $\Gamma \cup \{ \varphi \} \subseteq Fm(\vdash)$. Similarly, $\vdash$ is said to be \textit{almost inconsistent} if it lacks theorems and $\Gamma \vdash \varphi$ for every $\Gamma \cup \{ \varphi \} \subseteq Fm(\vdash)$ such that $\Gamma \ne \emptyset$. The following result is part of the folklore.

\begin{Lemma}\label{Lem:Inconsistent-models}
A logic $\vdash$ is inconsistent (resp.\ almost inconsistent) if and only if $\ModS(\vdash)$ is the class of isomorphic copies of $\langle {\bf 1}, \{ 1 \} \rangle$ (resp.\ of $\langle {\bf 1}, \{ 1 \} \rangle$ and $\langle {\bf 1}, \emptyset \rangle$).
\end{Lemma}

The lemma easily implies that any two inconsistent (resp.\ almost inconsistent) logics are equi-interpretable (since any translation between their languages is necessarily an interpretation).

\begin{Corollary}
The class of all inconsistent (resp.\ almost inconsistent) logics is a member of $\class{Log}$.
\end{Corollary}

In the light of the above corollary, the main results of this section can be summarized as follows:

\begin{Theorem}\label{Thm:extrema}
The poset $\class{Log}$ lacks a minimum.\ Moreover, its maximum is the class of all inconsistent logics, and its unique coatom is the class $\mathcal{K}$ of all almost inconsistent logics. In particular, a logic $\vdash$ lacks theorems if and only if $\llbracket \vdash \rrbracket \leq {\mathcal{K}}$.
\end{Theorem}

\begin{proof}
We first prove that $\class{Log}$ has no minimum. Suppose, with a view to contradiction, that $\class{Log}$ has a minimum  $\llbracket \vdash \rrbracket$. Then let $\kappa \coloneqq \vert Fm(\vdash) \vert$ and consider the language $\LL$ that consists in $k^{+}$ binary connectives $\{ \multimap_{\alpha} \colon \alpha <  \kappa^{+}\}$. For every $\alpha < \kappa^{+}$, let $\A_{\alpha}$ be the $\LL$-algebra with universe $\{ 1, 0, a \}$ and operations defined for every $p , q \in A$ and $\beta < \kappa^{+}$ as follows:
\begin{displaymath}
p \multimap_{\beta} q\coloneqq \left\{\begin{array}{@{\,}ll}
1 & \text{if $p=q$ or $\beta \ne \alpha$}\\
0 & \text{if $p \ne q$ and $\beta = \alpha$.}\\
\end{array} \right.
\end{displaymath}
Let also $\vdash_{\kappa^{+}}$ be the logic (formulated in a countable set of variables) induced by the class of \textit{reduced} matrices
\[
\class{M} \coloneqq \{ \langle \A_{\alpha}, \{ 1, a \}\rangle \colon \alpha < \kappa^{+} \}.
\]
Clearly,  $\class{M} \subseteq \ModS(\vdash_{\kappa^{+}})$.

Since $\vdash$ is the minimum of $\class{Log}$, there is an interpretation $\btau$ of $\vdash$ into $\vdash_{\kappa^{+}}$. On cardinality grounds, there is $\alpha < \kappa^{+}$ such that the symbol $\multimap_{\alpha}$ does not occur in the formulas $\{ \btau(\varphi) \colon \varphi \in \LL \}$. In particular, this implies that the matrix $\langle\A_{\alpha}^{\btau}, \{ 1, a \} \rangle$ is not reduced.

On the other hand, we know that $\tarski_{\vdash}^{\A_{\alpha}^{\btau}} \{ 1, a \}$ is the identity relation, since $\btau$ is an interpretation of $\vdash$ into $\vdash_{\kappa^{+}}$ and $\langle \A_{\alpha}, \{ 1, a \} \rangle \in \ModS(\vdash_{\kappa^{+}})$. Now, since $A_{\alpha}= \{ 1, 0, a \}$, the only deductive filter of $\vdash$ on $\A_{\alpha}^{\btau}$ extending properly $\{ 1, a \}$ is forcefully $\{ 1, 0, a \}$. Hence we obtain that
\[
\tarski_{\vdash}^{\A_{\alpha}^{\btau}} \{ 1, a \} = \leibniz^{\A_{\alpha}^{\btau}} \{ 1, a \} \cap \leibniz^{\A_{\alpha}^{\btau}} \{ 1, 0, a \} = \leibniz^{\A_{\alpha}^{\btau}} \{ 1, a \} \cap A_{\alpha}^{2} = \leibniz^{\A_{\alpha}^{\btau}} \{ 1, a \}.
\]
But this implies that $\leibniz^{\A_{\alpha}^{\btau}} \{ 1, a \}$ is the identity relation, which is false.

Now, from Lemma \ref{Lem:Inconsistent-models} it follows easily that the class of inconsistent (resp.\ almost inconsistent) logics is the maximum (resp.\ a coatom) of $\class{Log}$. Hence, in order to establish the second part of the theorem 
 it only remains to show that the class $\mathcal{K}$  of all almost inconsistent logics is the unique coatom of $\class{Log}$, and that a logic $\vdash$ lacks theorems if and only if $\llbracket \vdash \rrbracket \leq {\mathcal{K}}$.

By Lemmas \ref{Lem:no-theorems-trivial}(ii) and \ref{Lem:Inconsistent-models} a logic $\vdash$ lacks theorems if and only if $\llbracket \vdash \rrbracket \leq {\mathcal{K}}$.
Hence it only remains to show that $\mathcal{K}$ is the unique coatom of $\class{Log}$. Suppose, with a view to contradiction, that there is a coatom $\llbracket \vdash \rrbracket$ in $\class{Log}$ such that $\vdash$ is not almost inconsistent. Since $\llbracket \vdash \rrbracket$ is neither the maximum of $\class{Log}$, nor comparable with $\mathcal{K}$, we know that $\vdash$ is not inconsistent and that it has theorems. Then there is a matrix $\langle \A, F\rangle \in \ModS(\vdash)$ such that $F \in  \mathcal{P}(A) \smallsetminus \{ \emptyset, A \}$. In particular, this implies that $\vert A \vert \geq 2$. Now, consider the matrix
\[
\langle \B, G \rangle \coloneqq \langle \A, F \rangle^{\vert A \vert} \in \PPP(\ModS(\vdash)) \subseteq \ModS(\vdash)
\]
and observe that $\vert B \vert > \vert A \vert$ by Cantor's Theorem.

Let  $\B^{+}$ be the expansion of $\B$ with all finitary operations on $B$, and consider the logic $\vdash^{+}$ formulated in $\vert Fm(\vdash) \vert$ variables induced by the matrix $\langle \B^{+}, G \rangle$. 

Bering in mind that all finitary operations on $B$ are term-function of $\B^{+}$, it is not hard to see that the matrix $\langle \B^{+}, G \rangle$ is reduced and that the logic $\vdash^{+}$ is equivalential (see \cite[Lemma 3.2]{Mo15a} if necessary). Moreover, we have that $\SSS(\B^{+}) = \{ \B^{+} \}$. Together with Proposition \ref{Prop:alg-translations}, this implies that the identity map is a translation of $\vdash$ into $\vdash^{+}$. Since $\vdash$ is a coatom of $\class{Log}$, this implies that either $\vdash^{+}$ is inconsistent or it is equi-interpretable with $\vdash$. As $G \ne B$ and $\langle \B^{+}, G \rangle$ is a model of $\vdash^{+}$, we know that $\vdash^{+}$ is not inconsistent, whence ${\vdash^{+}} \leq {\vdash}$.

Together with the fact that $\langle \A, F \rangle \in \ModS(\vdash)$, this implies that $\ModS(\vdash^{+})$ contains a matrix of size $\vert A \vert$. However, from Lemma \ref{Lem:Equivalent_Alg_Sem} it follows that every non-trivial member of $\ModS(\vdash^{+})$ has cardinality $\geq \vert B \vert$. Together with the fact that $\vert A \vert < \vert B \vert$, this implies that $\A$ is trivial, which is false.
\end{proof}

\begin{Remark}
The proof above of the first part of Theorem \ref{Thm:extrema} suggests that the lack of a minimum in $\class{Log}$ can be amended if we impose restrictions on the cardinality of the languages in which logics are formulated.\footnote{The reader may have noticed that also the proof that finite suprema need not exist in $\class{Log}$ relies on the fact that the cardinality of languages in which logics are formulated is unbounded. However, in that case, it is not clear to the authors that imposing cardinality restriction on the size of the languages would be sufficient to recover the existence of suprema in $\class{Log}$.} To be more precise, we will show that the following poset has a minimum for every infinite cardinal $\kappa$: 
\[
\class{Log}_{\kappa} \coloneqq \{ \llbracket \vdash \rrbracket \colon \vert \LL_{\vdash} \vert \leq \kappa \} \subseteq \class{Log}.
\]

To this end, recall that the \textit{basic logic} $\vdash_{\class{V}}$ of a variety $\class{V}$ \cite{FMo14c,Mo13a} is the logic in the language of $\class{V}$ (formulated in a countable set of variables) induced by the following class of matrices
\[
\{ \langle \A, F\rangle \colon \A \in \class{V} \text{ and }F \subseteq A \}.
\]
Given an infinite cardinal $\kappa$, we consider the language $\LL_{\kappa}$ comprising $\kappa$ different $n$-ary symbols for every $n \in \omega$. Then let $\class{V}_{\kappa}$ be the variety of all $\LL_{\kappa}$-algebras. Clearly $\llbracket \vdash_{\class{V}_{\kappa}}\rrbracket \in \class{Log}_{\kappa}$. More interestingly, 
we shall prove that $\llbracket \vdash_{\class{V}_{\kappa}}\rrbracket$ is indeed the minimum of $\class{Log}_{\kappa}$.

Consider a logic $\vdash$ such that $\vert \LL_{\vdash} \vert \leq \kappa$. We can assume without loss of generality that the language of $\vdash$ is of size $\kappa$. Then there is a surjective translation $\btau \colon \LL_{\kappa} \to \LL_{\vdash}$. We will show that $\btau$ is an interpretation of $\vdash_{\class{V}_{\kappa}}$ into $\vdash$. To this end, consider $\langle \A, F \rangle \in \RRR(\Mod(\vdash))$. Since $\btau$ is surjective, the algebras $\A$ and $\A^{\btau}$ are term-equivalent. In particular, this implies that the matrix $\langle \A^{\btau}, F \rangle$ is reduced. Together with the fact that $\A^{\btau} \in \class{V}_{\kappa}$, this implies that $\langle \A^{\btau}, F \rangle \in \ModS(\vdash_{\class{V}_{\kappa}})$. Hence, with an application of Proposition \ref{Prop:interpretation2}, we conclude that $\btau$ is an interpretation.
\qed
\end{Remark}

As a consequence  of the remark we have: 
 
\begin{Corollary}
The upset of $\class{Log}$ generated by the set $\{ \llbracket \vdash_{\class{V}_{\kappa}}\rrbracket \colon \kappa \text{ is an infinite cardinal}\}$ is $\class{Log}$.
\end{Corollary}


\section{Relations with the lattice of varieties}

For $\kappa \in \omega$, a $k$\textit{-deductive system} $\vdash$ \cite{BP92} is a consequence relation over $Fm_{\LL}(\kappa)^{k}$ (for some language $\LL$ and infinite cardinal $\kappa$) that, moreover, is substitution invariant in the sense that for every substitution $\sigma$,
\[
\text{if }\Gamma \vdash \langle \varphi_{1}, \dots, \varphi_{k}\rangle \text{, then }\{ \langle \sigma(\gamma_{1}), \dots, \sigma(\gamma_{k})\rangle \colon \langle \gamma_{1}, \dots, \gamma_{k}\rangle \in \Gamma \} \vdash \langle \sigma(\varphi_{1}), \dots, \sigma(\varphi_{k}) \rangle
\]
for every $\Gamma \cup \{ \langle \varphi_{1}, \dots, \varphi_{k}\rangle \} \subseteq Fm_{\LL}(\kappa)^{k}$. 

\begin{exa}

Observe that $1$-deductive systems coincide with logics. Moreover, every variety $\class{K}$ can be associated with a $2$-deductive system $\vDash_{\class{K}}$ formulated over $Fm(\omega)^{2}$ as follows. For every $\Gamma \cup \langle \varphi, \psi \rangle \subseteq Fm(\omega)^{2}$ we set
\begin{align*}
\Gamma \vDash_{\class{K}} \langle \varphi, \psi \rangle \Longleftrightarrow& \text{ for every $\A \in \class{K}$ and homomorphism }f \colon \Fm(\omega) \to \A\\
& \text{ if }f(\epsilon) = f(\delta) \text{ for all }\langle \epsilon, \delta \rangle \in \Gamma\text{, then }f(\varphi) = f(\psi).
\end{align*}
The relation $\vDash_{\class{K}}$ is a notational variant of the standard equational consequence relative to $\class{K}$ (formulated in countably many variables).
\qed
\end{exa}

The theory of $k$-deductive systems is a smooth generalization of that of logics (for the details, see for instance \cite{BP92,Ra06}). In particular, every $k$-deductive system $\vdash$ can be associated with a class $\ModS(\vdash)$ of models of the form $\langle \A, F \rangle$ where $\A$ is an algebra and $F \subseteq A^{k}$. Bearing this in mind, we say that an \textit{interpretation} of a $k$-deductive system $\vdash$ into another $k$-deductive system $\vdash'$ is a translation $\btau$ of the language of $\vdash$ into that of $\vdash'$ such that $\langle \A^{\btau}, F \rangle \in \ModS(\vdash)$, for every $\langle \A, F \rangle \in \Mod(\vdash')$. We denote by $\class{Syst}(k)$ the poset of classes $\llbracket \vdash \rrbracket$ of equi-interpretable $k$-deductive systems, ordered under interpretability. Let also $\class{Equiv}(k)$ be the subposet of $\class{Syst}(k)$ that contains the classes $\llbracket \vdash \rrbracket$ such that $\vdash$ is an equivalential $k$-deductive system. A straightforward adaptation of the proof of Theorem \ref{Thm:ProtoalgebraicLattice} shows that $\class{Equiv}(k)$ is a set-complete lattice.

Recall form the Introduction that a variety $\class{K}$ is  \textit{interpretable} \cite{Tay73} into another variety $\class{V}$, when $\class{V}$ is term-equivalent to some variety $\class{V}^{\ast}$ whose reducts (in a smaller signature) belong to $\class{K}$, in which case  we write ${\class{K}} \leq {\class{V}}$. When ${\class{K}} \leq {\class{V}}$ and ${\class{V}} \leq {\class{K}}$ we say that ${\class{K}}$  and ${\class{V}}$ are \textit{equi-interpretable}. The class of all varieties equi-interpretable with $\class{K}$ is denoted by 
$\llbracket \class{K}\rrbracket $ and is called the \textit{interpretability type} of $\class{K}$. Moreover, we denote by $\class{Var}$ the lattice of interpretability types of varieties ordered by the relation $\leq$  defined as follows: $\llbracket \class{K}\rrbracket \leq \llbracket \class{V}\rrbracket $ if and only if ${\class{K}} \leq {\class{V}}$. The next result draws a relation between  $\class{Syst}(2)$, and the lattice $\class{Var}$ of interpretability types of varieties.

\begin{Proposition}
The map given by the rule $\llbracket \class{K}\rrbracket \longmapsto \llbracket \vDash_{\class{K}} \rrbracket$ is a lattice-embedding of $\class{Var}$ into $\class{Equiv}(2)$. 
\end{Proposition}

\begin{proof}[Proof sketch.]
It is well known that if $\class{K}$ is a variety, then $\vDash_{\class{K}}$ is an algebraizable \cite{BP89} (and, therefore, equivalential) $2$-deductive system such that
\[
\ModS(\vDash_{\class{K}}) = \{ \langle \A, \{ \langle a, a \rangle \colon a \in A \} \rangle \colon \A \in \class{K} \}.
\]
As a consequence, a variety $\class{K}$ is interpretable into another variety $\class{V}$ if and only if $\vDash_{\class{K}}$ is interpretable into $\vDash_{\class{V}}$. Hence, the map in the statement is an order-embedding. The fact that it is a lattice homomorphism follows from the description of infima and suprema in $\class{Var}$ \cite{GaTa84,Neum74}, and from a straightforward  adaptation of the description of infima and suprema of equivalential logics given here to the case of $2$-deductive systems.
\end{proof}

The above result gives a logical explanation of some known facts about $\class{Var}$. For instance, the fact that $\class{Var}$ is a lattice (as opposed to a poset only) can be viewed as a consequence of the fact that equivalential $2$-deductive systems form a lattice.\ Similarly, the fact that $\class{Var}$ has no coatoms \cite[Chpt.\ 2]{GaTa84} follows from a variant of Theorem \ref{Thm:extrema}, and the observation that every two-deductive system of the form $\vDash_{\class{K}}$ has at least one theorem, namely $\langle x, x \rangle$.

We conclude this section by showing that there is a meet-homomorphism from $\class{Var}$ to $\class{Log}$ (Theorem \ref{Thm:Mu_Map}). To this end, given a language $\mathcal{L}$ we denote by $T_{m}(\mathcal{L})$ the set of all its $m$-ary terms
in the variables $x_1, \ldots, x_m$. Then, for every $\mathcal{L}$-algebra $\A$ and $n>0$, the $n$\textit{-th matrix power} of $\A$ is the algebra
\[
\A^{[n]} \coloneqq \langle A^{n}; \{ m_{t} \colon t \in T_{kn}(\A)^{n} \text{ for some positive }k \in \omega \},
\]
where for each $t = \langle t_{1}, \dots, t_{n}\rangle \in T_{kn}(\A)^{n}$, we define $m_{t} \colon (A^{n})^{k} \to A^{n}$ as follows: if $a_{j} = \langle a_{j1}, \dots, a_{jn}\rangle \in A^{n}$ for $j = 1, \dots, k$, then
\[
m_{t}(a_{1}, \dots, a_{k}) = \langle t^{\A}_{i}(a_{11}, \dots, a_{1n}, \dots, a_{k1}, \dots, a_{kn}) :1 \leq  i \leq n \rangle.
\]

For $0<n\in\omega$, the $n$\textit{-th matrix power} of a class $\mathsf{K}$ of similar algebras is the class
$\mathsf{K}^{[n]} \coloneqq \III \{ \A^{[n]} : \A \in \mathsf{K} \}$. Applications of the matrix power construction range from the algebraic study of category equivalences and adjunctions \cite{DaWe83,Mc96,Mor16b} to the study of clones \cite{Ne70}, Maltsev conditions \cite{Ta75a,GaTa84}, and finite algebras \cite{HoMcKe88}. 

Given a variety $\class{K}$, we denote by $\vdash_{\class{K}}^{2}$ the logic formulated in countably many variables induced by the class of matrices
\begin{equation}\label{Eq:matrix-power}
\III \{ \langle \A^{[2]}, \{ \langle a, a \rangle \colon a \in A \}\rangle \colon \A \in \class{K}\}.
\end{equation}
We rely on the following observation \cite[Thm.\ 8]{MorRaf18-prevarieties}:

\begin{Theorem}\label{Thm:on-prevarieties}
If $\class{K}$ is variety, then $\ModS(\vdash^{2}_{\class{K}})$ is the class in (\ref{Eq:matrix-power}).
\end{Theorem}

As a consequence we obtain the desired result  (cf.\ \cite[Prop.\ 7]{GaTa84}):

\begin{Theorem}\label{Thm:Mu_Map}
The map defined by the rule $\llbracket \class{K}\rrbracket \longmapsto\llbracket \vdash_{\class{K}}^{[2]}\rrbracket$ is a meet-homomorphism from $\class{Var}$ into $\class{Log}$.
\end{Theorem}

\begin{proof}
We claim that if $\class{K}$ and $\class{V}$ are varieties such that ${\class{K}} \leq {\class{V}}$, then ${\vdash_{\class{K}}^{[2]}} \leq {\vdash_{\class{V}}^{[2]}}$. To prove this, let $\btau$ be an interpretation of $\class{K}$ into $\class{V}$. It is not hard to see that the map $\btau^{\ast}$, defined by the rule $\langle t_{1}, t_{2}\rangle \longmapsto \langle \btau(t_{1}), \btau(t_{2})\rangle$, is an interpretation of $\class{K}^{[2]}$ into $\class{V}^{[2]}$. We will show that $\btau^{\ast}$ is also an interpretation of $\vdash_{\class{K}}^{[2]}$ into $\vdash_{\class{V}}^{[2]}$. To this end, consider $\langle \A, F \rangle \in \ModS(\vdash_{\class{V}}^{[2]})$. By Theorem \ref{Thm:on-prevarieties} there is $\B \in \class{V}^{[2]}$ such that $\langle \A, F \rangle \cong \langle \B, \{ \langle b, b \rangle \colon b \in B \} \rangle$. As $\btau^{\ast}$ is an interpretation of $\class{K}^{[2]}$ and $\class{V}^{[2]}$, this yields $\B^{\btau^{\ast}} \in \class{K}^{[2]}$. Hence, by Theorem \ref{Thm:on-prevarieties} we obtain 
\[
\langle \A, F \rangle \in \III(\langle \B^{\btau^{\ast}}, \{ \langle b, b \rangle \colon b \in B \} \rangle) \subseteq \ModS(\vdash_{\class{K}}^{[2]}).
\]
We conclude that $\btau^{\ast}$ is an interpretation of $\vdash_{\class{K}}^{[2]}$ into $\vdash_{\class{V}}^{[2]}$, establishing the claim.

Let $\mu \colon \class{Var} \to \class{Log}$ be the map in the statement. From the claim it follows that $\mu$ is well-defined. Then we turn to prove that it is a meet-homomorphism. To this end, given two varieties $\class{K}$ and $\class{V}$, we set $\class{K} \bigotimes \class{V} \coloneqq \III \{ \A \bigotimes \B \colon \A \in \class{K} \text{ and }\B \in \class{V} \}$. It is easy to see that $\class{K} \bigotimes \class{V}$ is a variety. Moreover, recall that $\llbracket \class{K} \bigotimes \class{V} \rrbracket$ is the meet of $\llbracket \class{K}\rrbracket$ and $\llbracket\class{V} \rrbracket$ in $\class{Var}$, and that $(\class{K} \bigotimes \class{V})^{[2]}$ is term-equivalent to $\class{K}^{[2]} \bigotimes \class{V}^{[2]}$ (see for instance \cite{GaTa84}). Together with Corollary \ref{Cor:FiniteCase} and Theorem \ref{Thm:on-prevarieties}, this implies that $\vdash_{\class{K} \bigotimes \class{V}}^{[2]}$ is term-equivalent to $\vdash_{\class{K}}^{[2]}\bigotimes \vdash_{\class{V}}^{[2]}$ as well. Hence we have that
\begin{align*}
\mu(\llbracket \class{K}\rrbracket) \land \mu(\llbracket \class{V}\rrbracket) &= \llbracket 	\vdash_{\class{K}}^{[2]}\rrbracket \land \llbracket 	\vdash_{\class{V}}^{[2]}\rrbracket = \llbracket \vdash_{\class{K}}^{[2]}\bigotimes \vdash_{\class{V}}^{[2]} \rrbracket = \llbracket\vdash_{\class{K} \bigotimes \class{V}}^{[2]} \rrbracket \\
&= \mu (\llbracket \class{K} \bigotimes \class{V} \rrbracket) = \mu(\llbracket \class{K}\rrbracket \land \llbracket \class{V}\rrbracket).
\end{align*}
This shows that $\mu$ is a meet-homomorphism, as desired.
\end{proof}

\paragraph{\bfseries Acknowledgements.}
Thanks are due to James G.\ Raftery for rising the question about whether the theory of the Maltsev and Leibniz hierarchy could be, to some extent, unified. The second author was supported by  the grant CZ.$02$.$2$.$69$/$0$.$0$/$0$.$0$/$17$\_$050$/$0008361$, OPVVV M\v{S}MT, MSCA-IF Lidsk\'{e} zdroje v teoretick\'{e} informatice.

\section*{Appendix}

\subsection*{A}

Recall that $\vdash_{\lnot}$ is the negation fragment of classical propositional logic. The following result is part of the folklore:

\begin{Proposition}\label{Prop:Axioms-neg}
The logic $\vdash_{\lnot}$ is axiomatized by the following rules:
\[
x, \lnot x \rhd y \qquad x \rhd \lnot \lnot x \qquad \lnot \lnot x \rhd x.
\]
\end{Proposition}

\begin{Theorem}\label{Thm:Fragment}
$\ModS(\vdash_{\lnot})$ is the class of matrices $\langle \A, F\rangle$ such that either $\A$ is trivial or ($\A \in \class{NA}$ and either $F= \emptyset$ or $F= \{ a \}$ for some $a \in A$ that is not a fixed point of $\lnot$).
\end{Theorem}

\begin{proof}
We begin by proving the inclusion from left to right. To this end, observe that $\vdash_{\lnot}$ is determined the the matrix $\langle {\bf 2}, \{ 1 \} \rangle$ where ${\bf 2}$ is the negation reduct of the two-element Boolean algebra with universe $\{ 0, 1 \}$. Then consider a matrix $\langle \A, F \rangle \in \ModS(\vdash_{\lnot})$ such that $\A$ is non-trivial. First we show that $\A \in \class{NA}$. The fact that ${\bf 2} \vDash x \thickapprox \lnot \lnot x$, together with Corollary \ref{Cor:Mod-Suszko}, implies that $\A \vDash x \thickapprox \lnot \lnot x$.

 It only remains to prove that $\A$ has at most one fixed point of $\lnot$. Suppose that  $a, b \in A$ are fixed points of $\lnot$. We prove that $\textup{Fg}^{\A}_{\vdash_{\lnot}}(F, p(a)) = \textup{Fg}^{\A}_{\vdash_{\lnot}}(F, p(b))$ for every unary polynomial function $p$ of $\A$. This implies that $\langle a, b \rangle \in \tarski_{\vdash_{\lnot}}^{\A}F$ by Proposition \ref{Prop:Polynomial}(ii)  and, since  $\langle \A, F \rangle \in \ModS(\vdash_{\lnot})$, we that $a = b$. Let $p$ be a unary polynomial function of $\A$. Because of the poor language of $\A$, every polynomial function $q(x)$ of $\A$ has the form
\[
q(x) =\underbrace{\lnot \dots \lnot}_{n\text{-times}}x \quad \text{ or } \quad q(x) =\underbrace{\lnot \dots \lnot}_{n\text{-times}}c
\]
for some $n \in \omega$ and $c \in A$. If $p$ is of the first shape, then since $a, b$ are fixed points of $\lnot$, it easily follows that $p(a) = a$ and $p(b) = b$; then using  the rule $x, \lnot x \rhd y$, which holds in the logic, it follows that $b \in \textup{Fg}^{\A}_{\vdash_{\lnot}}(F, a)$ and $a \in \textup{Fg}^{\A}_{\vdash_{\lnot}}(F, b)$ and, therefore, $\textup{Fg}^{\A}_{\vdash_{\lnot}}(F, p(a)) = \textup{Fg}^{\A}_{\vdash_{\lnot}}(F, p(b))$. If $p$ is of the second shape, then there is nothing to prove. 
 Hence we conclude that $\A$ is a negation algebra.

Now, recall that $\vdash_{\lnot}$ is determined by a matrix $\langle {\bf 2}, \{ 1 \} \rangle$, whose set of designated elements is a singleton. By a  minor variant of \cite[Thm.\ 8]{AFRM15}, this implies that $\ModS(\vdash_{\lnot})$ is a class of matrices $\langle \A, F \rangle$ such that $F$ is either empty or a singleton. Then consider a matrix $\langle \A, F \rangle \in \ModS(\vdash_{\lnot})$ such that $\A$ is non-trivial. We know that $\A$ is a negation algebra and that $F$ is either empty or a singleton. Suppose, with a view to contradiction, that $F = \{ a \}$ for a fixed point $a$ of $\lnot$. Since $x, \lnot x \vdash_{\lnot} y$, this implies that $F = A$ and, therefore, that $\A$ is trivial which is false. This establishes the inclusion from left to right.

To prove the inclusion from right to left, consider a matrix $\langle \A, F \rangle$ in the right-hand side of the display in the statement. If $\A = {\bf 1}$, then either $F = \emptyset$ or $F = \{ 1 \}$. In both cases, $\langle \A, F \rangle \in \ModS(\vdash_{\lnot})$, since $\vdash_{\lnot}$ has no theorems. Then we suppose that $\A$ is non-trivial, in which case $\A \in \class{NA}$ and either $F = \emptyset$ or $F = \{ a \}$ for some $a \in A$ that is not a fixed point of $\lnot$. Together with Proposition \ref{Prop:Axioms-neg}, this implies that $\langle \A, F \rangle \in \Mod(\vdash_{\lnot})$. 

It only remains to prove that $\tarski_{\vdash_{\lnot}}^{\A}F$ is the identity relation. To prove this, consider two different elements $b, c \in A$. First we consider the case where $F = \emptyset$. Since $\A$ is a negation algebra, it has at most one fixed point of $\lnot$. Thus we can assume without loss of generality that $b$ is not a fixed point of $\lnot$. Together with Proposition \ref{Prop:Axioms-neg}, this implies that $\langle \A, \{ b \} \rangle \in \Mod(\vdash_{\lnot})$. Moreover, clearly we have that $F = \emptyset \subseteq \{ b \}$ and $\langle b, c \rangle \notin \leibniz^{\A}\{ b \} \subseteq \tarski_{\vdash_{\lnot}}^{\A} \emptyset = \tarski_{\vdash_{\lnot}}^{\A} F$. 

Then we consider the case in which $F = \{ a \}$ for some $a \in A$ that is not a fixed point of $\lnot$. Since $\A$ has at most one fixed point of $\lnot$, $\A \vDash x \thickapprox \lnot \lnot x$, and $b \ne c$ one of the following conditions holds:
\benroman
\item Either ($b \ne \lnot b$ and $b \ne \lnot a$) or ($c \ne \lnot c$ and $c \ne \lnot a$).
\item Either ($b = \lnot b$ and $c = \lnot a$) or ($c = \lnot c$ and $b = \lnot a$).
\eroman
If condition (i) holds, we can assume without loss of generality that $b \ne \lnot b$ and $b \ne \lnot a$. If $c = a$, then $\langle b, c \rangle \notin \leibniz^{\A}\{ a \} \subseteq \tarski_{\vdash_{\lnot}}^{\A} \{ a \} = \tarski_{\vdash_{\lnot}}^{\A} F$. Then consider the case where $c \ne a$. By Proposition \ref{Prop:Axioms-neg} we know that $\{ a, b \}$ is a deductive filter of $\vdash_{\lnot}$. Hence we have that $\langle b, c \rangle \notin \leibniz^{\A}\{ a, b \} \subseteq \tarski_{\vdash_{\lnot}}^{\A} \{ a \} = \tarski_{\vdash_{\lnot}}^{\A} F$. Then suppose that the condition (ii) holds. We can assume without loss of generality that $b = \lnot b$ and $c = \lnot a$. In this case we have that $\lnot b \notin \{ a \}$ and $\lnot c \in \{ a \}$ which, by Proposition \ref{Prop:Polynomial}(ii), implies that $\langle b, c \rangle \notin\tarski_{\vdash_{\lnot}}^{\A} F$.

This concludes the proof that $\tarski_{\vdash_{\lnot}}^{\A}F$ is the identity relation.
\end{proof}

\subsection*{B} 
We close the paper with an observation on languages with constant symbols.
If a logic $\vdash$ has constants in its language, then we can obtain a new language by keeping all the connectives of $\LL_\vdash$ and replacing each constant $c$ by a  unary operation $\ast_c$. Then we can transform every algebra $\A$ for the language $\LL_\vdash$ into an algebra $\A^{co}$ for the new language, where $\ast_c^{\A^{co}}$ is the unary constant map to $c^{\A}$. The logic $\vdash^{co}$ in the new language induced by the class of matrices  
\[
\{\langle \A^{co}, F\rangle:  \langle \A, F\rangle \in \Mod(\vdash)\}
\]
 is the incarnation of the logic $\vdash$ in our  setting of logics with languages without constants. 
 
Note that if the language of $\vdash$ has no constant symbols, then ${\vdash} = {\vdash^{co}}$. It is therefore natural to say of any two logics $\vdash$ and $\vdash'$, possibly with constants,  that $\vdash$ is \textit{interpretable} into $\vdash'$ if $\vdash^{co}$ is interpretable into $\vdash'^{co}$ in the sense of Definition \ref{def:interpretation}. 
Alternatively, and with similar ideas to the ones used in the paper,  the reader can easily figure out how to modify our notion of a concrete interpretation to accommodate interpretations between languages possibly with constants.
%
%
%

\bibliographystyle{plain}

\end{document}